\documentclass[12pt]{article}
\usepackage[a3paper]{geometry}
\usepackage[numbers]{natbib}
\usepackage{graphicx}
\usepackage{color, multirow,amsfonts,amsmath}
\usepackage{url}
\usepackage{rotating}
\begin{document}
\title{Mixed Poisson process with Stacy mixing variable - Working version
\bigskip
\\ Pavlina K. Jordanova
\\{\small{\it Faculty of Mathematics and Informatics, Konstantin Preslavsky University of Shumen, 115 "Universitetska" str., 9712 Shumen, Bulgaria, Corresponding author: pavlina\_kj@abv.bg}}
\\ Mladen Savov,
\\{\small{Faculty of Mathematics and Informatics, Sofia University, and Institute of Mathematics and Informatics, Bulgarian Academy of Sciences, Sofia, Bulgaria}}
\\ Assen Tchorbadjieff
\\{\small{Institute of Mathematics and Informatics, Bulgarian Academy of Sciences, Sofia, Bulgaria}}
\\ Milan Stehl\'ik
\\{\small{Institute of Applied Statistics and Linz Institute of Technology, Johannes Kepler University, Linz, Austria and
\\ Institute of Statistics, University of Valpara\'iso, and Facultad de Ingenier\'{\i}a, Universidad Andr\'es Bello, Valpara\'iso, Chile}}}
\maketitle

\begin{abstract}
Stacy distribution defined for the first time in 1961 provides a flexible framework for modelling of a wide range of real-life behaviours. It appears under different names in the scientific literature and contains many useful particular cases. Homogeneous Poisson processes are appropriate apriori models for the number of renewals up to a given time $t>0$. This paper mixes them and considers a Mixed Poisson process with Stacy mixing variable. We call it a Poisson-Stacy process. The resulting counting process is one of the Generalised Negative Binomial processes, and the distribution of its time-intersections are very-well investigated in the scientific literature. Here we define and investigate their joint probability distributions. Then, the corresponding mixed renewal process is investigated and Exp-Stacy and Erlang-Stacy distributions are defined and partially studied.

The paper finishes with a simulation study of these stochastic processes. Some plots of the probability density functions, probability mass functions, mean square regressions and sample paths are drawn together with the corresponding code for the symulations.
\end{abstract}


\section{DESCRIPTION OF THE MODEL AND PRELIMINARIES}

In 1920 the authors of Greenwood and Yule \cite{GreenwoodYule} introduce a mixed Poisson process with Gamma mixing variable. In 2016 \cite{LMJ} Jordanova and Stehlik mix a homogeneous Poisson process with Pareto mixing variable and investigate its properties. Although that model covers a large variety of real data situations it has too high overdispersion. Here we consider a similar model with moderate overdispersion. More precisely the Pareto distribution in Jordanova and Stehlik \cite{LMJ} is replaced with Stacy distribution. The last one seems to be defined in 1961 by Stacy \cite{Stacy}  via its probability density function (p.d.f.)
\begin{equation}\label{pdfStacy}
P_\xi(x) = |c|\frac{\beta^{c\alpha}}{\Gamma(\alpha)}x^{c\alpha-1}e^{-(\beta x)^c}, \quad x \geq 0.
\end{equation}
Here and hereafter we suppose that $c \not = 0$, $\alpha> 0$, and denote by $\Gamma(\alpha) := \int_0^{\infty} y^{\alpha - 1} e^{-y} dy$. This distribution provides a flexible modelling of a wide range of real life behaviours and has different names. For example, it captures the Three-parameter Gamma, Generalized Gamma, Amoroso distribution (without shifting), among others. Stoyanov \cite{JStoyanov} calls this distribution hyper-exponential. He considers Berg's problem \cite{Berg}, and shows that, under different parametrization, when $c \geq \frac{1}{2}$, this distribution is determined via its moments, and when $c \in \left(0, \frac{1}{2}\right)$ it is moment indeterminate. There are many different Generalized Gamma distributions, therefore, we prefer the name of the author of (\ref{pdfStacy}), and the following parametrization $\xi \in Stacy(\alpha, \beta; c)$. Stacy and Mihram \cite{Stacy1965} showed that this distribution is closed with respect to power transformations and more precisely, for $\delta \in \mathbb{R}$, $\xi^\delta \in Stacy(\alpha, \beta^\delta; \frac{c}{\delta})$. This distribution possesses the following scaling property for $\lambda > 0$, $\lambda\xi \in Stacy(\alpha, \frac{\beta}{\lambda}; c)$. When $\lambda < 0$ it is switched to the negative half semi axis. This class of probability laws contains many useful and popular particular cases. It is obvious that $Stacy(\alpha, \beta; 1)$ coincides with $Gamma (\alpha, \beta)$ distribution. More generally,
\begin{equation}\label{StacyGamma}
\xi \in Stacy(\alpha, \beta; c) \iff \xi \stackrel{d}{=} \eta^{1/c}, \quad \eta \in Gamma (\alpha, \beta^c).
\end{equation}
The notation
$\eta \in Gamma (\alpha, \beta^c)$ means that
$$P_{\eta}(x) = \frac{\beta^{\alpha c}}{\Gamma(\alpha)}x^{\alpha-1}e^{-\beta^c x}, \quad x > 0.$$

For $\alpha = 1$ and
\begin{description}
  \item[$\checkmark$] $c > 0$ this distribution is Positive Weibull with shape parameter equal to $c$ and scale parameter equal to $\beta$.
  \item[$\checkmark$] $c < 0$ this distribution is Fr$\acute{e}$chet with shape parameter equal to $-c$ and scale parameter equal to $\beta$.
\end{description}

The exponential, Erlang, chi-square, Weibull, folded normal, Fr$\acute{e}$chet, their powers and inverse versions, among others are also Stacy distributed. For more particular cases see for example Stacy and Mihram \cite{Stacy1965}, Johnson et al. \cite{Johnsonetal}, Dadpay et al. \cite{Dadpay}, Crooks \cite{Crooks} or Khodabin and Ahmadabadi \cite{Khodabin}. The third author investigates the information properties of this family. Stacy and Mihram \cite{Stacy1965} found that for all $\frac{\delta}{c} > -\alpha$,
\begin{equation}\label{moments}
\mathbb{E}(\xi^\delta) = \frac{\Gamma\left(\alpha + \frac{\delta}{c}\right)}{\beta^\delta\Gamma(\alpha)}
\end{equation}
and otherwise $\mathbb{E}\xi^\delta = \infty$. The corresponding cumulative distribution
function (c.d.f.) is
\begin{equation}\label{cdf}
F_\xi(x):= \mathbb{P}(\xi < x) = \left\{\begin{array}{ccc}
                                       \frac{\gamma((\beta x)^c, \alpha)}{\Gamma(\alpha)} & , & x \geq 0, \,\, c > 0 \\
                                        1 - \frac{\gamma((\beta x)^c, \alpha)}{\Gamma(\alpha)} & , & x \geq 0, \,\, c < 0.
                                      \end{array}
\right.
\end{equation}

The exact likelihood ratio tests for homogeneity and the scale parameter of this distribution have been derived in \cite{Stehlik2008}.

 Here and hereafter, we denote by $\gamma(x, \alpha) :=\int_0^{x} y^{\alpha - 1} e^{-y} dy$ the lower incomplete Gamma function, and by
\begin{equation}\label{Kratzel}
Z_\rho^\nu(x):= \int_0^\infty y^{\nu-1}e^{-(y^\rho+\frac{x}{y})}dy = \int_0^\infty t^{-\nu-1}e^{-(t^{-\rho}+tx)}dt, \quad x > 0,
\end{equation}
$$\rho \in (0, \infty),\,\, \nu \in \mathbb{C},\,\, {\text{and}} \,\,{\text{when}}\,\,\rho\leq 0, \,\,\mathbb{R}e(\nu) < 0,
$$
the {\bf  Kr$\ddot{a}$tzel function} investigated for example in Kr$\ddot{a}$tzel \cite{Kratzel}, Kilbas et al. \cite{Kilbasetal} and Princy \cite{Princy2014}. The notation $\eta \in Exp(1)$ means that the random variable (r.v.) $\eta$ is exponentially distributed with $\mathbb{E}\eta = 1$. $P_\theta(x)$ stands for the p.d.f. of the r.v. $\theta$. Note that $Z_\rho^\nu(0) = \frac{1}{|\rho|}\Gamma\left(\frac{\nu}{\rho}\right)$ when $\rho\nu > 0$.

It is straightforward to see that if $\mathbb{R}e(\nu) < 0$, then
\begin{description}
  \item[$\checkmark$] for $x > 0$, $Z_{0}^\nu(x) = \frac{x^\nu}{e} \Gamma(-\nu)$.
  \item[$\checkmark$] For $x > -1$, $Z_{-1}^\nu(x) = (1+x)^\nu \Gamma(-\nu)$.
\end{description}

The Laplace-Stieltjes transform of Stacy distribution is,
\begin{equation}\label{StacyMGF}
\mathbb{E}(e^{-t\xi}) = |c|\frac{\beta^{c\alpha}}{\Gamma(\alpha)} \int_0^\infty x^{c\alpha-1}e^{-(\beta x)^c-tx}dx =  \frac{|c|}{ \Gamma(\alpha)}Z_{-c}^{-c\alpha}\left(\frac{t}{\beta}\right), \quad t > 0.
\end{equation}

The general theory of Mixed Poisson processes could be found for example in Grandel \cite{Grandel} or Karlis and Xekalaki \cite{KarlisandXekalaki}. Here the homogeneous Poisson process is mixed with a Stacy distribution and it is called {\bf Poisson-Stacy process}. It is one of the Generalised Negative Binomial processes, and the distribution of its time-intersections are very-well investigated in the scientific literature. Korolev et al. \cite{KorolevProbstat} call them {\bf GG-mixed Poisson distribution} and Korolev and Zeifman \cite{Kolorlevetal2019} call them {\bf Generalized negative binomial (GNB)}, prove that it is {\bf Mixed geometric} and describe the mixing r.v. when $\alpha \in (0, 1]$.  Here we define and investigate their joint probability distributions and distributions of their additive increments. The corresponding mixed renewal process is investigated and {\bf Exp-Stacy} and {\bf Erlang-Stacy distributions} are  partially studied.

The computational option for further implementation of the obtained results does not require significant quantity of additional resources. There are several different implementations of Stacy distribution in R software \cite{R1}. In following numerical experiments it is used Stacy distribution implementation in library "VGAM" \cite{R2}.The direct computations of results probability density functions (p.d.fs.) or mean square regression can be obtained directly after using standard numerical integration functionality for solving eq. 5 in R.
\bigskip

\section{EXP-STACY AND ERLANG-STACY DISTRIBUTIONS}

{\bf Definition 1.} \label{Def:1} We say that the r.v. $\tau$ is {\bf{Exp-Stacy distributed with parameters $\alpha > 0$, $\beta > 0$, and $c \in \mathbb{R}\backslash\{0\}$}}, if it has a p.d.f.
\begin{equation}\label{ExpStacyDensity}
P_\tau(t) = \left\{ \begin{array}{ccc}
                                   0 & , & t \leq 0\\
                                   \frac{|c|}{\beta\Gamma(\alpha)}Z_{-c}^{-c\alpha-1}\left(\frac{t}{\beta}\right)& , & t > 0
                      \end{array}\right..
\end{equation}

\begin{figure}
\includegraphics[scale=.75]{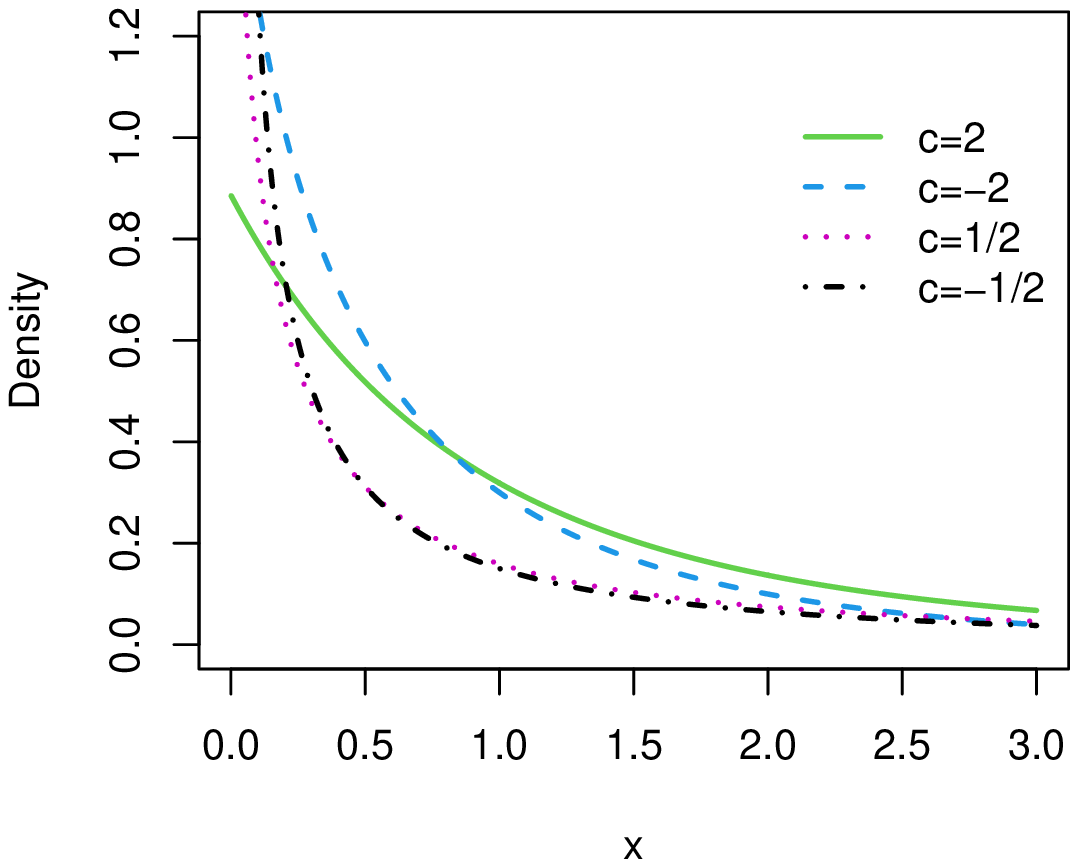} \hfill \includegraphics[scale=.75]{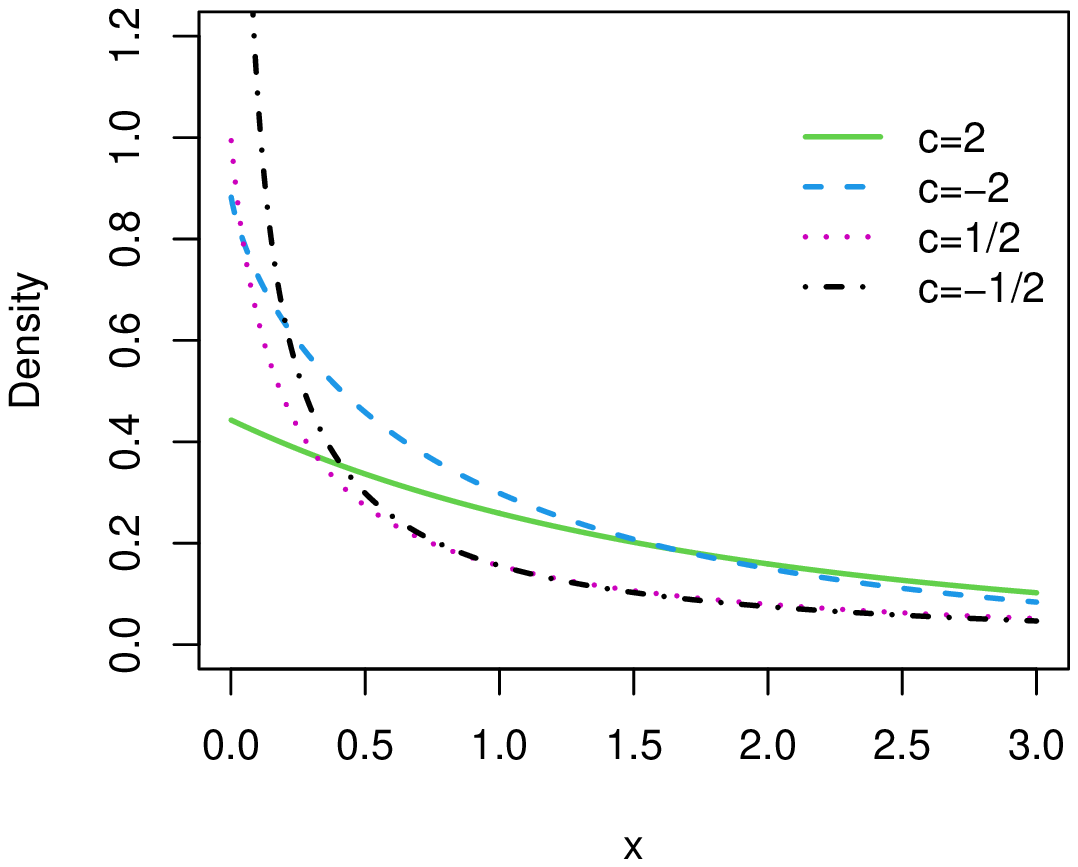}

\caption{Probability density functions of $Exp-St(1, 1, c = 2, -2, 0.5, -0.5)$(left) and $Exp-St(1, 2, c)$(right) distribution.\label{fig:ESt1}}
\end{figure}

\begin{figure}
\includegraphics[scale=.75]{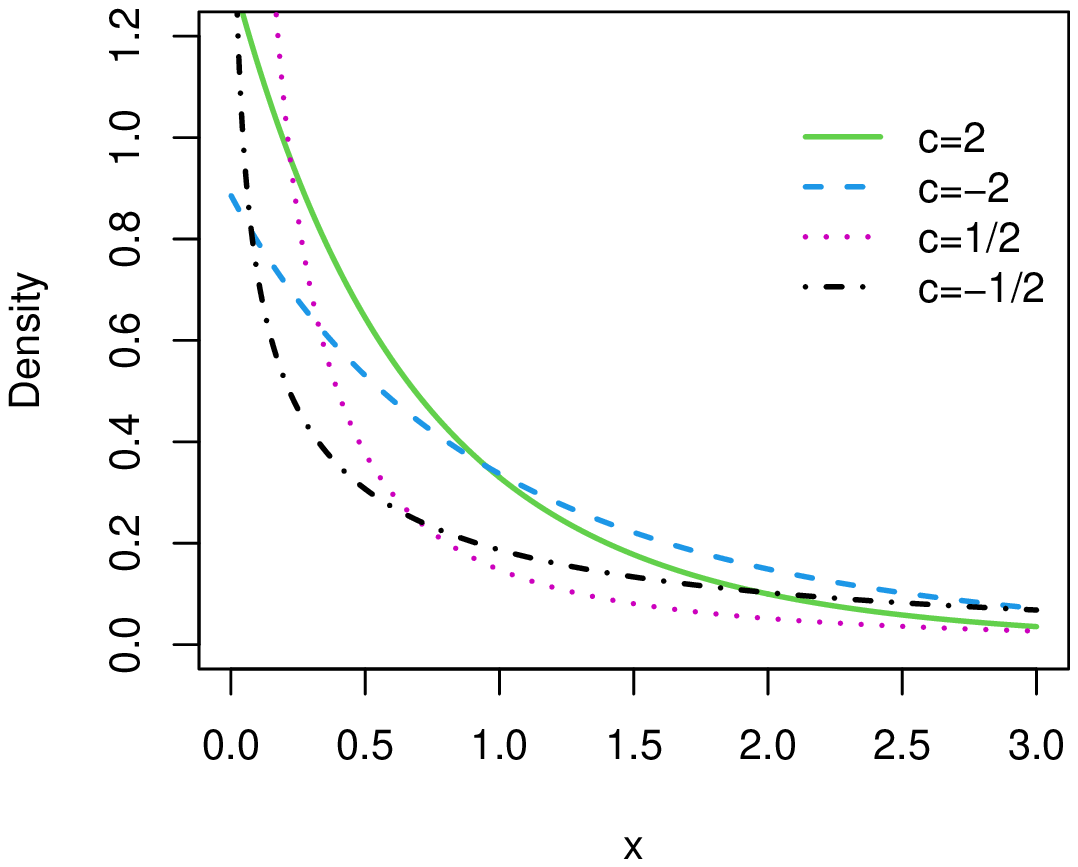} \hfill \includegraphics[scale=.75]{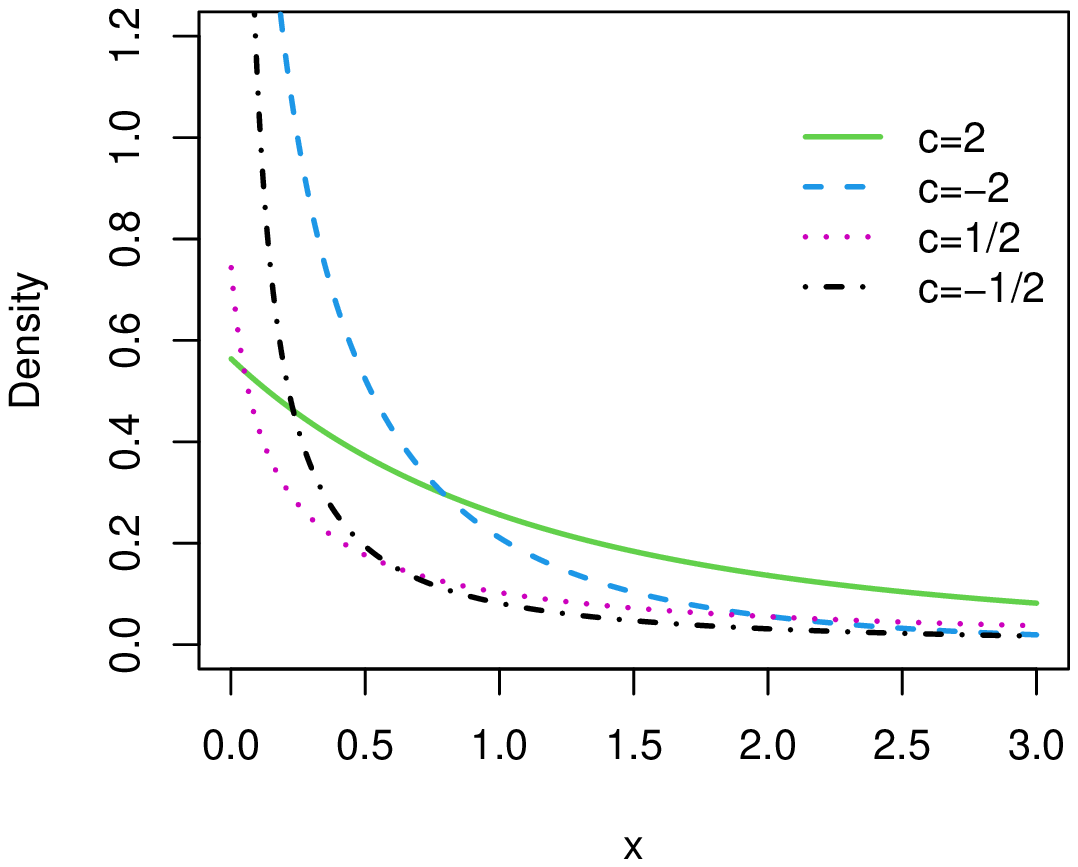}

\caption{Probability density functions of $Exp-St(2, 1, c)$ (left) and $Exp-St(1/2, 1, c)$(right) distribution.\label{fig:ESt2}}
\end{figure}

For brevity we denote this as  $\tau \in Exp-St(\alpha, \beta; c)$.
Following this notation the plots of p.d.f. for different parameters are computed and shown in Figure \ref{fig:ESt1} and Figure \ref{fig:ESt2}.
Note that in one of the demonstrated results, when $c = 1$,  $\tau$ is Pareto distributed due to the well-known fact that $P_\tau(t) = \frac{\alpha \beta^\alpha}{(t+\beta)^{\alpha + 1}}$ when $t > 0$, and $P_\tau(t) = 0$, otherwise.

On Figures \ref{fig:ESt1} and \ref{fig:ESt2} some plots of this probability density function for different values of its parameters are presented.

Another interesting property is that when $c = -1$, as it is discussed in Willmot \cite{Willmoth1993}. Then, $P_\tau(t) = 2\frac{\beta}{\Gamma(\alpha)}\sqrt{(\beta t)^{\alpha - 1}}K_{\alpha - 1}(2\sqrt{\beta t})$ when $t > 0$ and $P_\tau(t) = 0$, otherwise.

Finally we remark that as far as $\int_0^\infty Z_{-c}^{-c\alpha-1}\left(\frac{t}{\beta}\right)dt = \frac{\beta \Gamma(\alpha)}{|c|}$ the Exp-Stacy distribution is proper.

{\bf Definition 2.} \label{Def:1} We say that the random vector (rv.) $(\tau, \xi)$ has {\bf{bivariate Exp-Stacy distribution of $I^{-st}$ kind with parameters $\alpha > 0$, $\beta > 0$, and $c \in \mathbb{R}\backslash\{0\}$}}, if it has a joint p.d.f.
\begin{equation}\label{BivariateExpStacyDensity}
P_{\tau, \xi}(t,x) = \left\{ \begin{array}{ccc}
|c|t\frac{\beta^{c\alpha}}{\Gamma(\alpha)}x^{c\alpha-1}e^{-(\beta x)^c-tx}& , & t > 0, x > 0\\
0 & , & otherwise
                      \end{array}\right..
\end{equation}
Briefly we will denote this in this way  $(\tau, \xi) \in Exp-Stacy-I(\alpha, \beta; c)$.

Let us note that in 2014 Princy \cite{Princy2014} introduces  {\bf Kr$\ddot{a}$tzel distribution} and describes its relation with the Generalized Gamma distribution. In 2019 Kudryavtsev \cite{Kudryavtsev} defines {\bf Gamma-exponential distribution (GE distribution)}. Both these $l$-types\footnote{If a r.v. $\xi$ has some probability law, we call $l$-type the set of distributions of all r.vs. $a\xi + b$, where $a\in \mathbb{R}$ and $b \in \mathbb{R}$. $l$ comes from "linear".} are different from the Exp-Stacy one.

In the next theorem, which contains the first main result of our paper, we investigate the properties of this distribution.

{\bf Theorem 1.} For $\alpha > 0$, $\beta > 0$, and $c \in \mathbb{R}\backslash\{0\}$, if $\xi \in Stacy(\alpha, \beta; c)$ and for $\lambda > 0$, $(\tau|\xi=\lambda) \in Exp(\lambda)$, then:

\begin{description}
  \item[a)] $\tau \in Exp-St(\alpha, \beta; c)$;
  \item[b)]  If $k > 0$ is a constant, then $k\tau \in Exp-St(\alpha, k\beta; c)$;
  \item[c)] $\tau \stackrel{d}{=} \frac{\eta}{\xi}$, where $\eta \in Exp(1)$, and $\xi$ and $\eta$ are independent.
  \item[d)] For all $z \geq -1$, and $\frac{z}{c}\leq \alpha$, $\mathbb{E}(\tau^z) = \frac{\Gamma(z+1)\beta^z}{\Gamma(\alpha)}\Gamma\left(\alpha-\frac{z}{c}\right)$, and for $\max\left(\frac{2}{c},\frac{1}{c}\right)\leq \alpha$,
      $$\mathbb{D}\tau = \frac{\beta^2}{\Gamma(\alpha)}\left(2\Gamma\left(\alpha-\frac{2}{c}\right) - \frac{1}{\Gamma(\alpha)}\Gamma^2\left(\alpha-\frac{1}{c}\right)\right).$$
  \item[e)] The joint distribution of $\tau$ and $\xi$ is $(\tau, \xi) \in Exp-Stacy-I (\alpha, \beta; c)$ and $(\tau, \xi) \stackrel{d}{=} \left(\frac{\eta}{\xi}, \xi\right)$, where $\eta \in Exp(1)$, and $\xi$ and $\eta$ are independent.
  \item[f)] For all $t > 0$,
  $$P_{\xi}(x|\tau=t) = x^{c\alpha}e^{-tx-(\beta x)^c}\frac{\beta^{c\alpha+1}}{Z_{-c}^{-c\alpha-1}\left(\frac{t}{\beta}\right)}, \quad x > 0,$$
and $P_{\xi}(x|\tau=t) = 0$ when $x \leq 0$.
  \item[d)] The mean square regression $\mathbb{E}(\tau|\xi=x) = \frac{1}{x}$, $x > 0$.
   \item[h)] The mean square regression is $$\mathbb{E}(\xi|\tau=t) = \frac{Z_{-c}^{-c\alpha-2}\left(\frac{t}{\beta}\right)}{\beta Z_{-c}^{-c\alpha-1}\left(\frac{t}{\beta}\right)}, \quad t > 0.$$
\end{description}

\begin{figure}
\includegraphics[scale=.65]{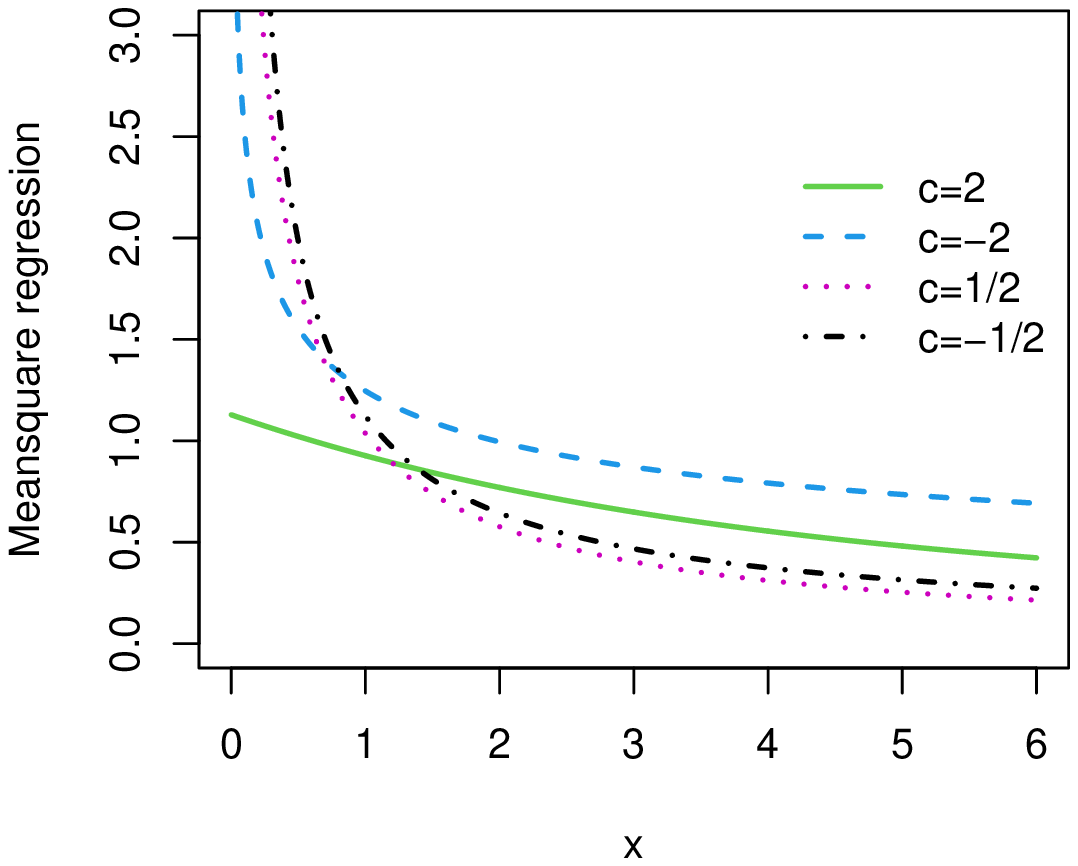} \hfill \includegraphics[scale=.65]{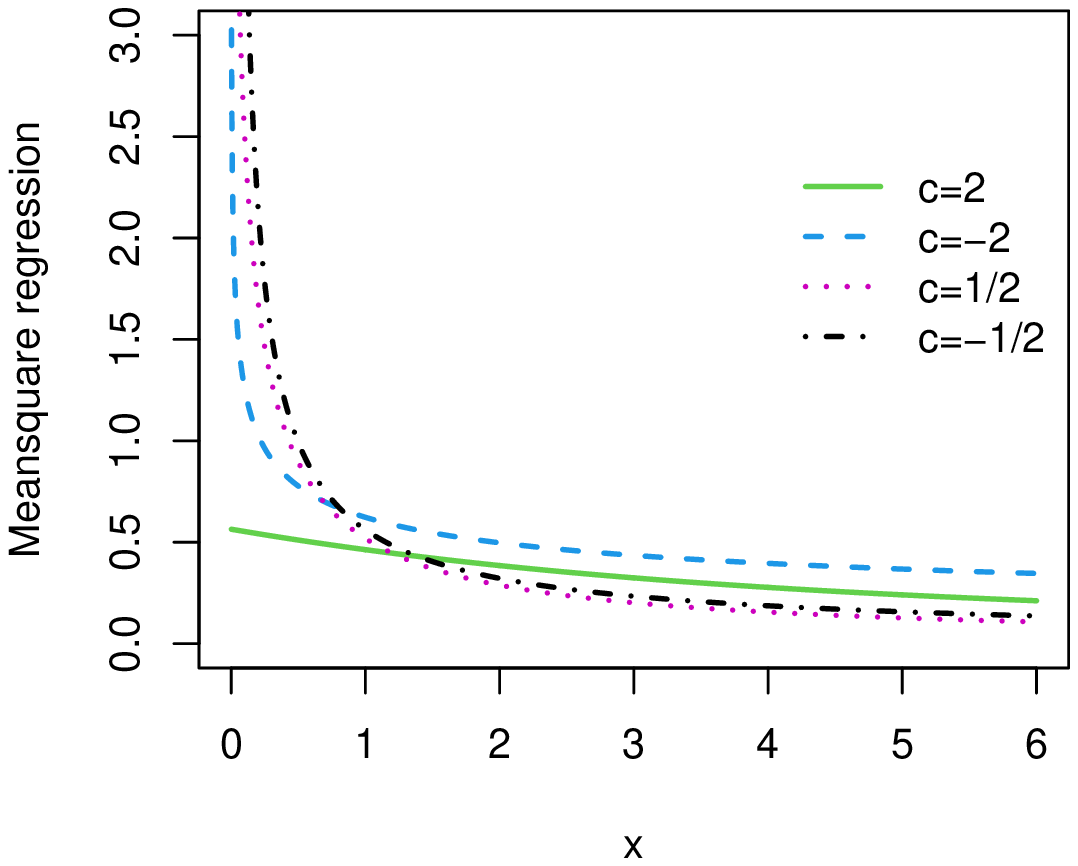}

\caption{Mean square regressions $\mathbb{E}(\xi|\tau=t)$, $t \geq 0$ of $Exp-St(1, 1, c = 2, -2, 0.5, -0.5)$(left) and $Exp-St(1, 2, c)$(right) distribution.\label{fig:MSRegressionESt1}}
\end{figure}

\begin{figure}
\includegraphics[scale=.75]{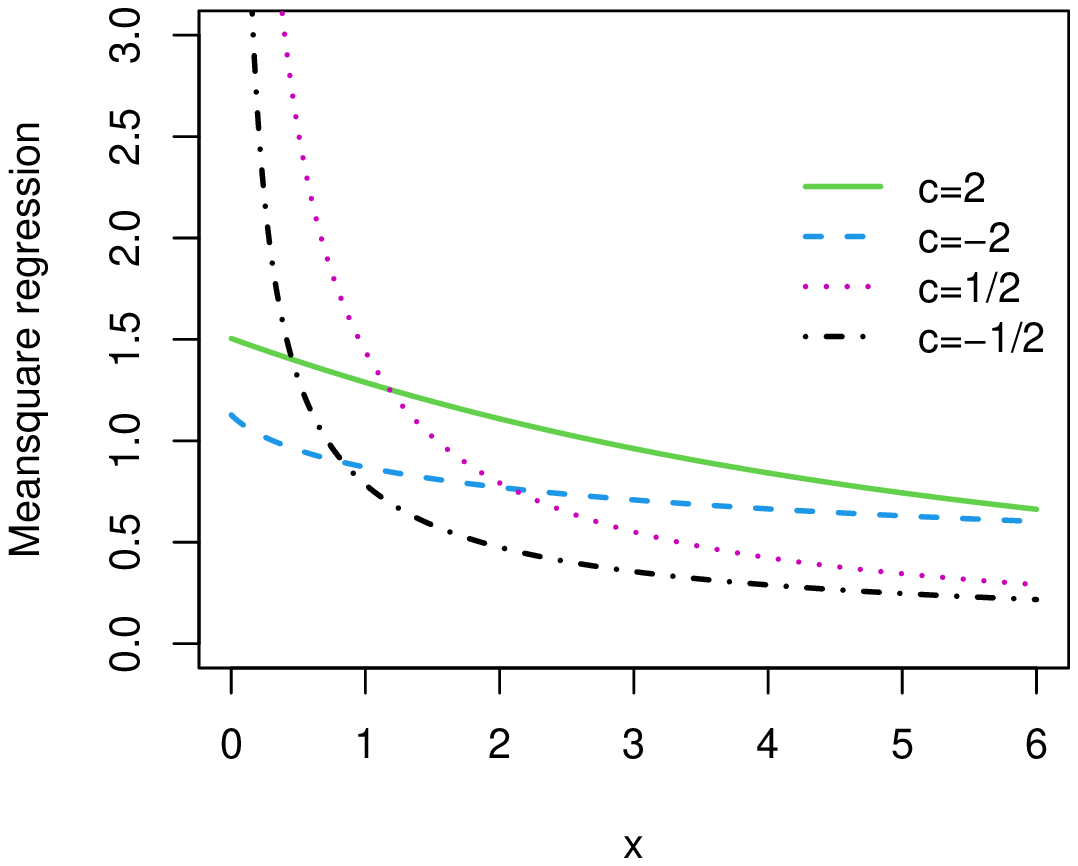} \hfill \includegraphics[scale=.75]{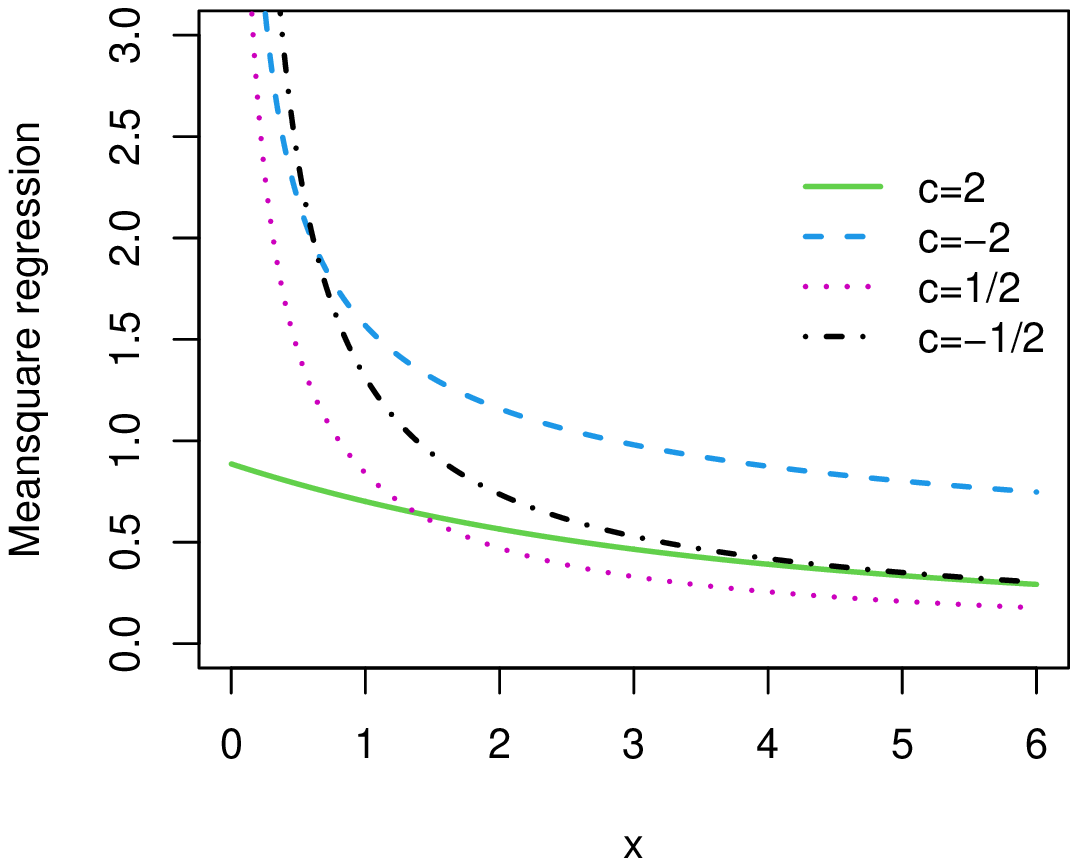}

\caption{Mean square regressions $\mathbb{E}(\xi|\tau=t)$, $t \geq 0$  of $Exp-St(2, 1, c)$ (left) and $Exp-St(1/2, 1, c)$(right) distribution.\label{fig:MSRegressionESt2}}
\end{figure}

{\bf Proof:} a) For $t > 0$ when substituting $y : = \beta x$ we obtain
\begin{eqnarray*}
  P_\tau(t) &=& \int_0^\infty P_\tau(t|\xi=x)P_\xi(x)dx = \int_0^\infty xe^{-xt}|c|\frac{\beta^{c\alpha}}{\Gamma(\alpha)}x^{c\alpha-1}e^{-(\beta x)^c}dx\\
   &=& |c|\frac{\beta^{c\alpha}}{\Gamma(\alpha)} \int_0^\infty x^{c\alpha}e^{-tx-(\beta x)^c}dx = \frac{|c|}{\beta\Gamma(\alpha)} \int_0^\infty y^{c\alpha}e^{-\frac{ty}{\beta}-y^c}dy = \frac{|c|}{\beta\Gamma(\alpha)}Z_{-c}^{-c\alpha-1}\left(\frac{t}{\beta}\right).
\end{eqnarray*}
Now, we compare it with (\ref{ExpStacyDensity}) and complete the proof of a).

b) For $k > 0$ and $x > 0$,
$$ P_{k\xi}(x) = \frac{1}{k}P_{\xi}\left(\frac{x}{k}\right) = \left\{ \begin{array}{ccc}
                                   0 & , & x \leq 0\\
                                   \frac{|c|}{k\beta\Gamma(\alpha)}Z_{-c}^{-c\alpha-1}\left(\frac{x}{k\beta}\right)& , & x > 0
                      \end{array}\right..$$
The rest follows by the uniqueness of the correspondence between p.d.f. and the distribution and formula (\ref{ExpStacyDensity}).

c) For $t > 0$, when substituting $y : = \beta x$ we obtain
\begin{eqnarray*}
  P_{\frac{\eta}{\xi}}(t) &=& \int_0^\infty P_{\frac{\eta}{\xi}}(t|\xi=x)P_\xi(x)dx = \int_0^\infty P_\frac{\eta}{x}(t)P_\xi(x)dx \\
  &=& \int_0^\infty xP_\eta(tx)P_\xi(x)dx =\int_0^\infty xe^{-xt}|c|\frac{\beta^{c\alpha}}{\Gamma(\alpha)}x^{c\alpha-1}e^{-(\beta x)^c}dx = \frac{|c|}{\beta\Gamma(\alpha)}Z_{-c}^{-c\alpha-1}\left(\frac{t}{\beta}\right).
\end{eqnarray*}
Now, we compare it with (\ref{ExpStacyDensity}) and by using a) complete the proof of a).

d) Consider $z \geq -1$ and $\frac{z}{c} \leq \alpha$. By consecutive application of the double expectation formula, $(\tau|\xi=\lambda) \in Exp(\lambda)$, the formula for the moments of the exponential distribution,  and (\ref{moments}) we get
$$\mathbb{E}(\tau^z) = \mathbb{E}(\mathbb{E}(\tau^z|\xi)) = \mathbb{E}\left(\frac{\Gamma(z+1)}{\xi^z}\right) = \Gamma(z+1)\mathbb{E}(\xi^{-z}) = \frac{\beta^z\Gamma(z+1)}{\Gamma(\alpha)}\Gamma\left(\alpha - \frac{z}{c}\right).$$

e) follows by the formula $P_{\tau, \xi}(t,x) =  P_{\tau}(t|\xi = x)P_{\xi}(x)$, when we replace the p.d.f. of the exponential distribution, use its scaling property, and formula (\ref{pdfStacy}).

f) By Bayes' formula for the densities, (\ref{pdfStacy}) and (\ref{ExpStacyDensity}), for $x > 0$ and $t > 0$,

$$P_{\xi}(x|\tau=t) = \frac{P_{\tau}(t|\xi = x)P_\xi(x)}{p_\tau(t)} = \frac{xe^{-tx}|c|\frac{\beta^{c\alpha}}{\Gamma(\alpha)}x^{c\alpha-1}e^{-(\beta x)^c}}{\frac{|c|}{\beta\Gamma(\alpha)}Z_{-c}^{-c\alpha-1}\left(\frac{t}{\beta}\right)} = x^{c\alpha}e^{-tx-(\beta x)^c}\frac{\beta^{c\alpha+1}}{Z_{-c}^{-c\alpha-1}\left(\frac{t}{\beta}\right)}$$
and $P_{\xi}(x|\tau=t) = 0$ when $x \leq 0$.

g) follows by the expectation of the Exponential distribution.

h) follows by the formula for the expectation, e) and (\ref{Kratzel}).
 \hfill $\Box$

{\it{Notes:}} For $c = 1$, the conclusion in f), means that for all $t > 0$, $(\xi|\tau = t) \in Gamma(\alpha+1, t+\beta)$, together with h) the last means that $\mathbb{E}(\xi|\tau=t) = \frac{\alpha+1}{t+\beta}$.

 {\bf Definition 3.} \label{Def:3} We say that the rv. $(\tau_1, \tau_2, ..., \tau_k)$ has {\bf{ Multivatiate Exp-Stacy distribution of $II^{-nd}$ kind with parameters $\alpha > 0$, $\beta > 0$, and $c \in \mathbb{R}\backslash\{0\}$}}, if it has a joint p.d.f.
\begin{equation}\label{MultivariateExpStacyDensityII}
P_{\tau_1, \tau_2, \ldots, \tau_k}(t_1, t_2, \ldots, t_k) =\left\{ \begin{array}{ccc}
                                                                      \frac{|c|}{\beta^k\Gamma(\alpha)}Z_{-c}^{-c\alpha-k}\left(\frac{t_1+t_2+\ldots+t_k}{\beta}\right)& , & t_1 > 0, t_2 > 0, \ldots, t_k > 0\\
                                                                      0 & , & otherwise
                      \end{array}\right..
\end{equation}
Briefly we will denote this in this way  $(\tau_1, \tau_2, \ldots, \tau_k) \in Exp-Stacy-II(\alpha, \beta; c)$.

In particular, when $c = 1$ we obtain
\begin{equation}\label{MultivariateExpStacyDensityIIc=1}
P_{\tau_1, \tau_2, \ldots, \tau_k}(t_1, t_2, \ldots, t_k) =\left\{ \begin{array}{ccc}
\frac{(\alpha+k-1)(\alpha+k-2)\ldots \alpha \beta^\alpha}{(\beta + t_1 + t_2 + \ldots + t_k)^{\alpha+k}}& , & t_1 > 0, t_2 > 0, \ldots, t_k > 0\\
                                                                      0 & , & otherwise
                      \end{array}\right..
\end{equation}

{\bf Definition 4.} \label{Def:4} We say that the r.v. $T_n$ is {\bf{Erlang-Stacy distributed with parameters $n \in \mathbb{N}$, $\alpha > 0$, $\beta > 0$, and $c \in \mathbb{R}\backslash\{0\}$}}, if it has a p.d.f.
\begin{equation}\label{ErlStacyDensity}
P_{T_n}(t) = \left\{ \begin{array}{ccc}
                                   0 & , & t \leq 0\\
                                   \frac{t^{n-1}|c|}{(n-1)!\beta^n\Gamma(\alpha)}Z_{-c}^{-c\alpha-n}\left(\frac{t}{\beta}\right)& , & t > 0
                      \end{array}\right..
\end{equation}

\begin{figure}
\includegraphics[scale=.75]{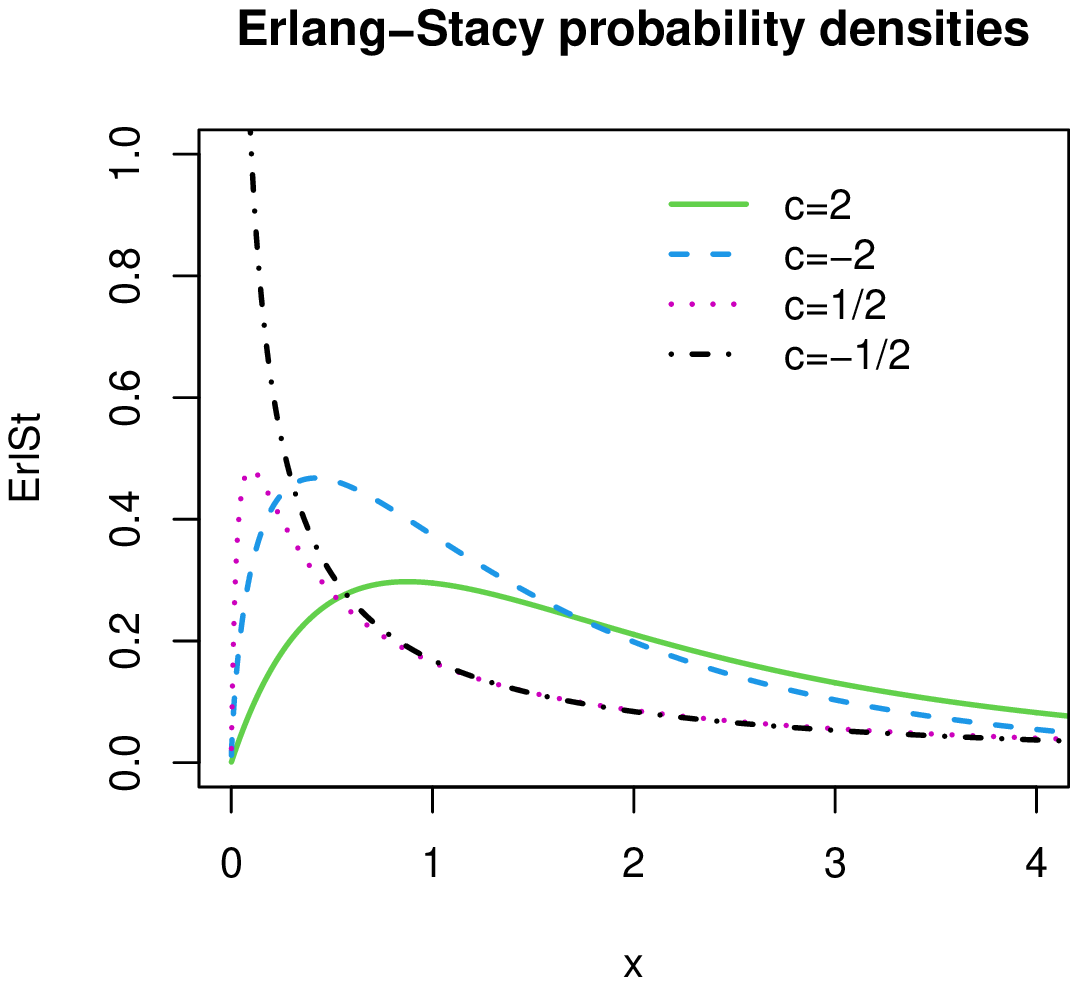} \hfill \includegraphics[scale=.75]{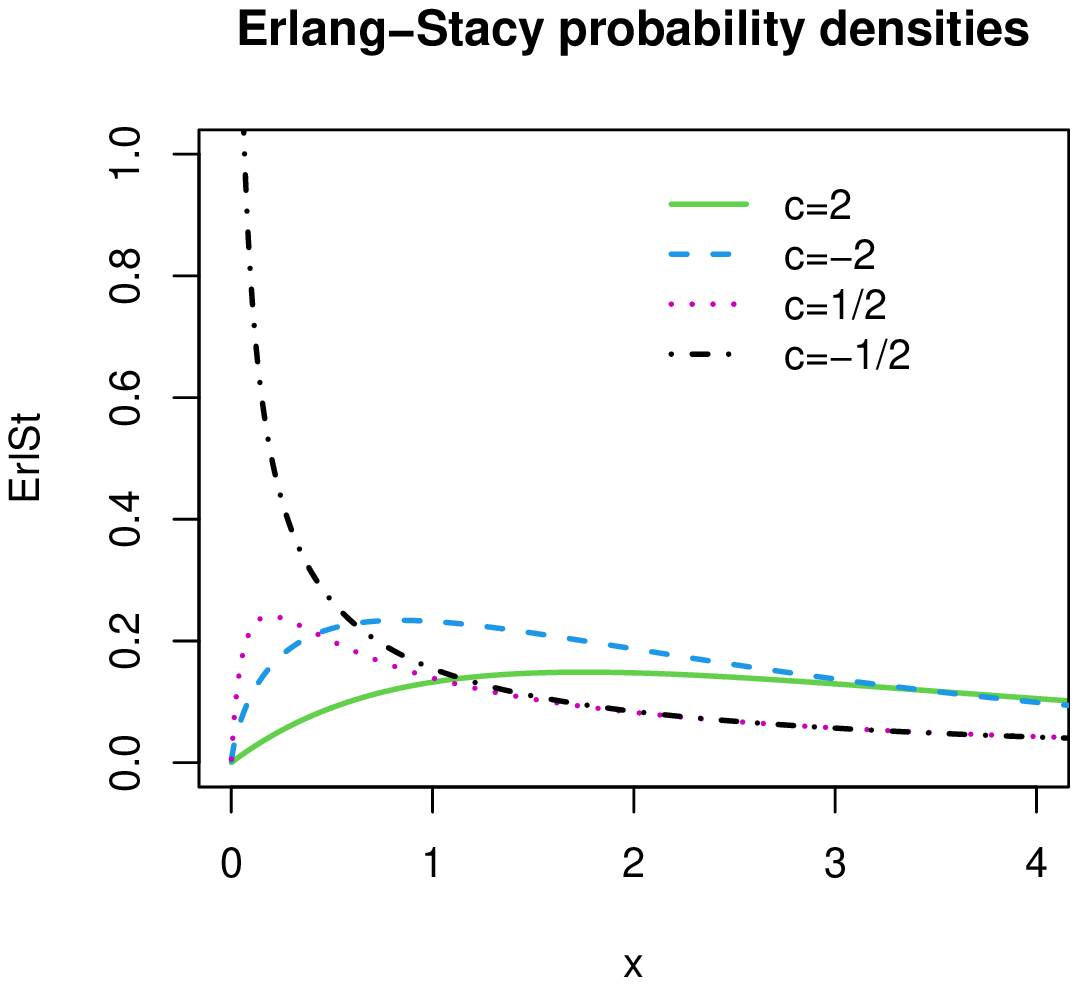}

\caption{Probability density functions of $Erlang-St(2; 1, 1, c)$(left) and $Erlang-St(2; 1, 2, c)$(right) distribution.\label{fig:ErlSt1}}
\end{figure}

\begin{figure}
\includegraphics[scale=.75]{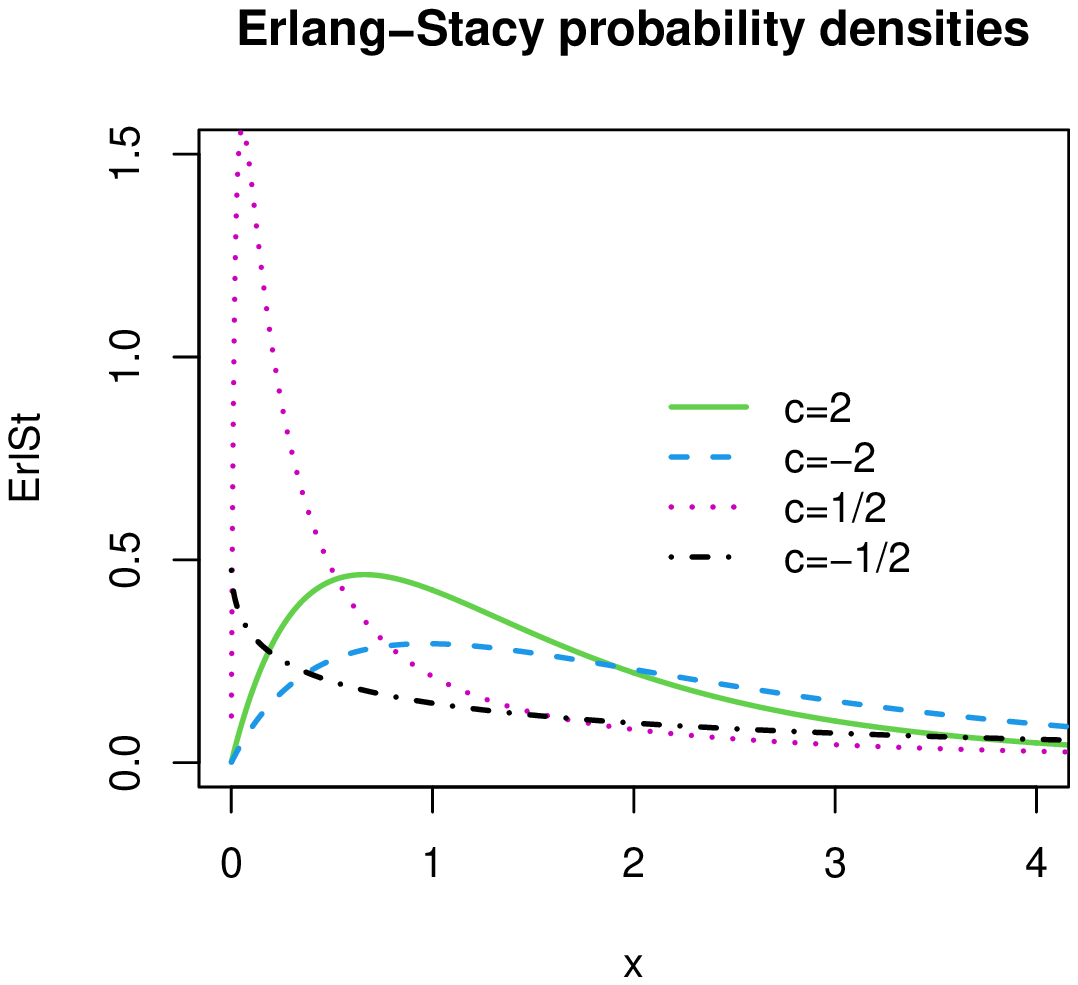} \hfill \includegraphics[scale=.75]{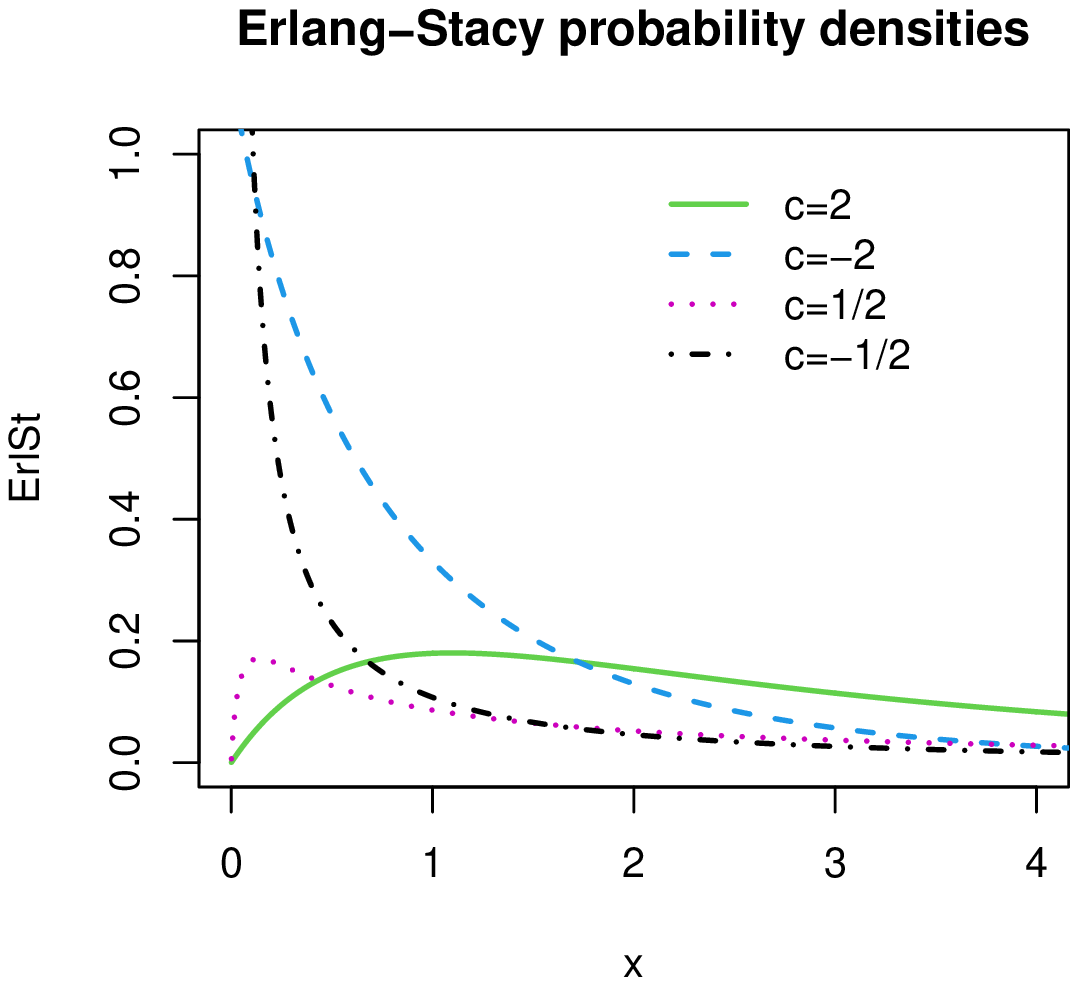}

\caption{Probability density functions of $Erlang-St(2; 2, 1, c)$ (left) and $Erlang-St(2; 1/2, 1, c)$(right) distribution.\label{fig:ErlSt2}}
\end{figure}

Briefly we will denote this in this way  $T_n \in Erlang-St(n;\alpha, \beta; c)$, $n \in \mathbb{N}$, $\alpha > 0$, $\beta > 0$, and $c \in \mathbb{R}\backslash\{0\}$. The p.d.f. for this distribution is plotted for different parameter values on Figures \ref{fig:ErlSt1} and \ref{fig:ErlSt2}.

In the next theorem, which is the second main result of this paper, we investigate some of the properties of Multivatiate Exp-Stacy distribution of $II^{-nd}$ kind and Erlang-Stacy distribution.

{\bf Theorem 2.} For $\alpha > 0$, $\beta > 0$, and $c \in \mathbb{R}\backslash\{0\}$, if $\xi \in Stacy(\alpha, \beta; c)$ and for $\lambda > 0$, $(\tau_1, \tau_2, ..., \tau_k|\xi=\lambda)$ are independent identically $Exp(\lambda)$ distributed r.vs., then,

\begin{description}
  \item[a)] $(\tau_1, \tau_2, \ldots, \tau_k) \in Exp-Stacy-II(\alpha, \beta; c)$.
  \item[b)] For $i = 1, 2, ..., k$, $\tau_i \in  Exp-Stacy(\alpha, \beta; c)$.
  \item[c)]$(\tau_1, \tau_2, \ldots, \tau_k)  \stackrel{d}{=} \left(\frac{\eta_1}{\xi}, \frac{\eta_2}{\xi}, \ldots, \frac{\eta_k}{\xi}\right)$, where $\eta_1, \eta_2, \ldots, \eta_k$ are i.i.d. $Exp(1)$, and independent on $\xi$.
 \item[d)] $T_n := \tau_1 + \ldots + \tau_n \in Erlang-St(n; \alpha, \beta; c)$,
 $${\text{for}}\,\,\frac{1}{c}\leq \alpha,\quad \mathbb{E}T_n = \frac{n\beta}{\Gamma(\alpha)}\Gamma\left(\alpha-\frac{1}{c}\right),$$
  $${\text{and}}\,\,{\text{for}}\,\,\max\left(\frac{1}{c}, \frac{2}{c}\right)\leq \alpha,\quad \mathbb{D}T_n = \frac{n\beta^2}{\Gamma(\alpha)}\left((n+1)\Gamma\left(\alpha-\frac{2}{c}\right) - \frac{n}{\Gamma(\alpha)}\Gamma^2\left(\alpha-\frac{1}{c}\right)\right).$$
 \item[e)] $T_n := \tau_1 + \ldots + \tau_n \stackrel{d}{=} \frac{\eta_1 + \eta_2 + \ldots + \eta_n}{\xi}$, where $\eta_1, \eta_2, \ldots, \eta_n$ are i.i.d. $Exp(1)$, and independent on $\xi$.
 \item[f)] $T_n \stackrel{d}{=} \frac{\theta_n}{\xi}$, where $\theta_n \in Gamma(n, 1)$ is independent on $\xi$.
  \item[g)]  If $k > 0$ is a constant, then $kT_n \in Erlang-St(n; \alpha, k\beta; c)$.
  \end{description}

{\bf Proof:} a) For $t_1, t_2, \ldots, t_k > 0$  by substituting $y : = \beta \lambda$ we obtain
\begin{eqnarray*}
P_{\tau_1, \tau_2, \ldots, \tau_k}(t_1, t_2, \ldots, t_k) &=& \int_0^\infty P_{\tau_1, \tau_2, \ldots, \tau_k}(t_1, t_2, \ldots, t_k|\xi = \lambda)P_{\xi}(\lambda)d\lambda = \int_0^\infty \lambda^k e^{-\lambda(t_1 + t_2 + \ldots + t_k)} |c|\frac{\beta^{c\alpha}}{\Gamma(\alpha)}\lambda^{c\alpha-1}e^{-(\beta \lambda)^c} d\lambda\\
&=& |c|\frac{\beta^{c\alpha}}{\Gamma(\alpha)} \int_0^\infty \lambda^{k+c\alpha-1} e^{-\lambda(t_1 + t_2 + \ldots + t_k) -(\beta \lambda)^c} d\lambda  = \frac{|c|}{\beta^k\Gamma(\alpha)} \int_0^\infty y^{k+c\alpha-1} e^{-y\frac{t_1 + t_2 + \ldots + t_k}{\beta} -y^c} dy\\
&=& \frac{|c|}{\beta^k\Gamma(\alpha)} Z_{-c}^{-c\alpha - k} \left(\frac{t_1 + t_2 + \ldots + t_k}{\beta}\right).
\end{eqnarray*}
In the case when some of $t_1, t_2, \ldots, t_k$ is negative the corresponding c.d.f. is 0. Therefore, $P_{\tau_1, \tau_2, \ldots, \tau_k}(t_1, t_2, \ldots, t_k) = 0$. Now, by comparing our expression to (\ref{MultivariateExpStacyDensityII}) we complete the proof.

b) follows by Theorem 1, a).

c) Analogously for $t_1, t_2, \ldots, t_k > 0$,
\begin{eqnarray*}
P_{\frac{\eta_1}{\xi}, \frac{\eta_2}{\xi}, \ldots, \frac{\eta_k}{\xi}}(t_1, t_2, \ldots, t_k) &=& \int_0^\infty P_{\frac{\eta_1}{\lambda}, \frac{\eta_2}{\lambda}, \ldots, \frac{\eta_k}{\lambda}}(t_1, t_2, \ldots, t_k|\xi = \lambda)P_{\xi}(\lambda)d\lambda \\
&=& \int_0^\infty P_{\eta_1, \eta_2, \ldots, \eta_k}(\lambda t_1, \lambda t_2, \ldots, \lambda t_k) \lambda^k P_{\xi}(\lambda)d\lambda \\
&=& \int_0^\infty \lambda^k e^{-\lambda(t_1 + t_2 + \ldots + t_k)} |c|\frac{\beta^{c\alpha}}{\Gamma(\alpha)}\lambda^{c\alpha-1}e^{-(\beta \lambda)^c} d\lambda = \frac{|c|}{\beta^k\Gamma(\alpha)} Z_{-c}^{-c\alpha - k} \left(\frac{t_1 + t_2 + \ldots + t_k}{\beta}\right).
\end{eqnarray*}
In the case when some of $t_1, t_2, \ldots, t_k$ is negative the corresponding c.d.f. is 0. Therefore, $P_{\tau_1, \tau_2, \ldots, \tau_k}(t_1, t_2, \ldots, t_k) = 0$. Now, by comparing the last expression to(\ref{MultivariateExpStacyDensityII}) we complete the proof.

d) $T_n = \tau_1 + \ldots + \tau_n \in Erlang-St(n; \alpha, \beta; c)$, follows by the integral form of the Total probability formula and the relation between the Erlang and Exponential distribution. The formula for $\mathbb{E}T_n$ is a consequence of Theorem 1, d) and the additive property of the expectation. The second moment of $T_n$ is obtained by the Double expectation formula. $\mathbb{E}(T_n^2) = \frac{n(n+1)}{\Gamma(\alpha)}\beta^2 \Gamma\left(\alpha - \frac{2}{c}\right)$. The well-known formula for the variance completes the proof.

e) is a corollary of c).

f) is a corollary of e) and the well-known presentation of the Erlang distribution as a convolution of Exponentials.

 b) For $k > 0$ and $x > 0$,
$$ P_{kT_n}(x) = \frac{1}{k}P_{T_n}\left(\frac{x}{k}\right) = \left\{ \begin{array}{ccc}
                                   0 & , & x \leq 0\\
                                   \frac{x^{n-1}|c|}{k^{n-1}(n-1)!k \beta^n\Gamma(\alpha)}Z_{-c}^{-c\alpha-n}\left(\frac{x}{k\beta}\right)& , & x > 0
                      \end{array}\right..$$
The rest follows by the uniqueness of the correspondence between p.d.f. and the distribution and formula (\ref{ErlStacyDensity}).
 \hfill $\Box$

{\it{Notes:}} For $c = 1$, $\frac{T_n}{\beta}$ has Beta prime distribution with scale parameter $n \in \mathbb{N}$ and shape parameter $\alpha > 0$,
\begin{equation*}\label{ExpStacyDensityc=1}
P_{T_n}(t) = \left\{ \begin{array}{ccc}
                                   0 & , & t \leq 0\\
                                   \frac{t^{n-1}\beta^\alpha}{(t+\beta)^{\alpha + n}{\text{B}}(n, \alpha)}& , & t > 0
                      \end{array}\right.; \quad \mathbb{E}T_n = \frac{n\beta}{\alpha - 1}, \alpha > 1; \quad \mathbb{D}T_n = \frac{n\beta^2(n+\alpha-1)}{(\alpha - 1)^2(\alpha-2)}, \alpha > 2,
\end{equation*}
where ${\text{B}}(n, \alpha) = \int_0^1 x^{n-1}(1-x)^{\alpha-1}dx$. The last means that this Beta prime distribution with integer scale parameter is a mixed Erlang distribution with Gamma mixing variable.

\bigskip
\section{THE MIXED POISSON-STACY PROCESS}

This section considers Mixed Poisson process (MPP) with Stacy mixing r.v. We call it {\bf Poisson-Stacy process}. It is time-changed {\bf Generalised Negative Binomial (GNB) process}. As far as it is a particular case of a MPP it is overdispersed. We concentrate our study mainly on the joint probability distributions, and on the distributions of the additive increments. The distribution of its time-intersections is very-well investigated in the scientific literature. Korolev et al. \cite{KorolevProbstat} call it {\bf GG-mixed Poisson distribution}, prove that when $\alpha, c \in (0, 1]$ it is a {\bf Mixed geometric} and describe the mixing r.v. The authors of Kudryavtsev \cite{Kudryavtsev} and Korolev and Zeifman \cite{Kolorlevetal2019} consider the case $c > 0$, and call the same distribution {\bf GNB distribution}.

\begin{figure}
\includegraphics[scale=.75]{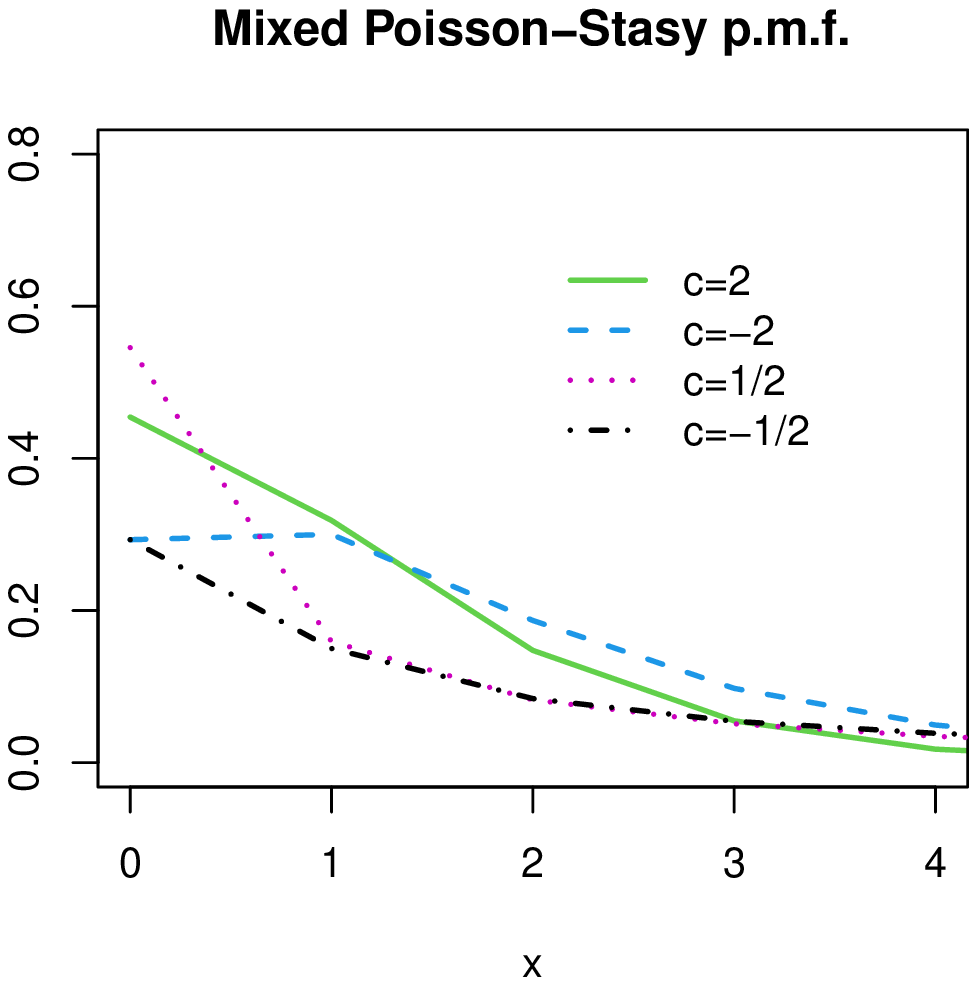} \hfill \includegraphics[scale=.75]{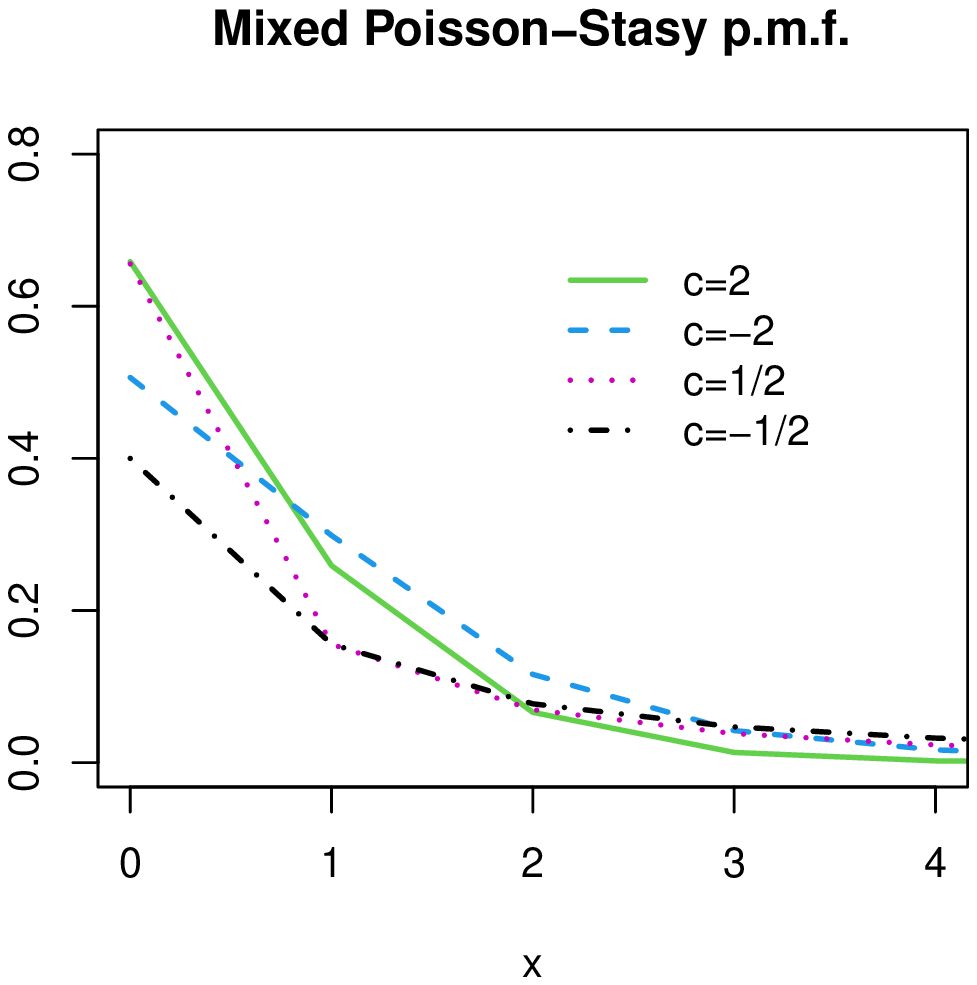}

\caption{Probability mass function of $MPStacy(1, 1, c)$(left) and $MPStacy(1, 2, c)$(right) distribution.\label{fig:MPSt1}}
\end{figure}

\begin{figure}
\includegraphics[scale=.75]{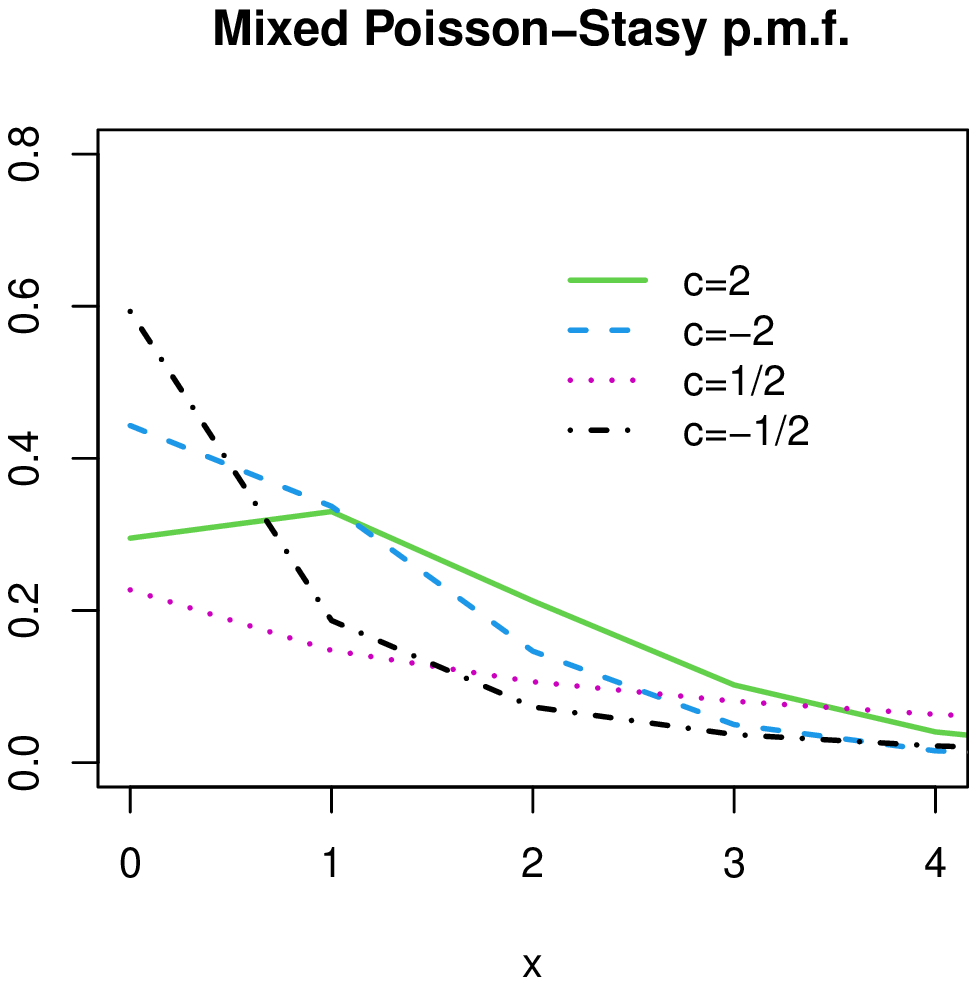} \hfill \includegraphics[scale=.75]{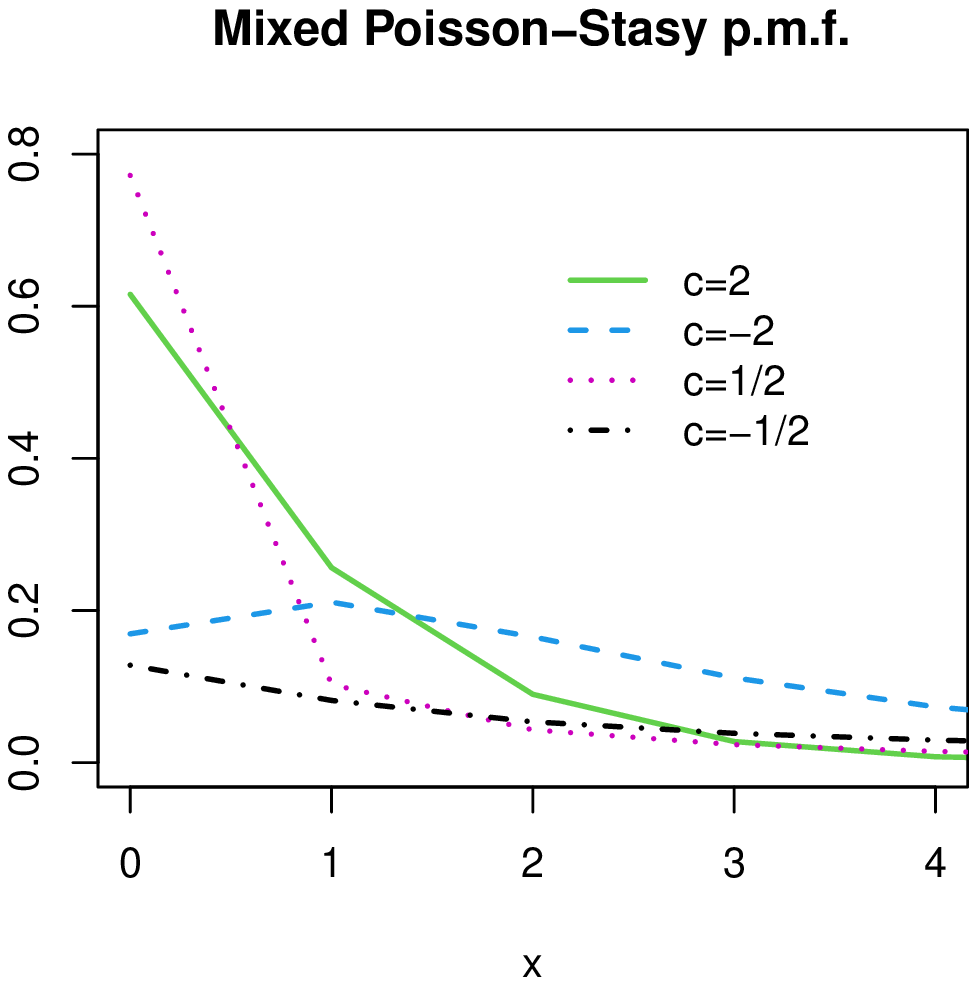}

\caption{Probability mass function of $MPStacy(2, 1, c)$ (left) and $MPStacy(1/2, 1, c)$(right) distribution.\label{fig:MPSt2}}
\end{figure}

{\bf Definition 5.} \label{Def:5} A r.v. $\theta$ has a {\bf{Mixed Poisson-Stacy distribution with parameters $\alpha > 0$, $\beta > 0$ and $c \in \mathbb{R}\backslash\{0\}$}} if
\begin{equation}\label{MixedPoissonStacy}
\mathbb{P}(\theta = n) = \frac{|c|}{n! \beta^n  \Gamma(\alpha)}Z_{-c}^{-c\alpha - n}\left(\frac{1}{\beta}\right), \quad n = 0, 1, \ldots.
\end{equation}
Briefly, $\theta \in MPStacy(\alpha, \beta, c)$. The p.m.fs. for different values of the parameters are shown in  Figures \ref{fig:MPSt1} and \ref{fig:MPSt2}.

{\bf Definition 6.} \label{Def:6} Let $\lambda(t): [0, \infty) \to [0, \infty)$ be a nonnegative, strictly increasing and continuous function, $\lambda(0) = 0$, $\xi \in Stacy(\alpha, \beta; c)$ and $N_1$ be a Homogeneous Poisson process with intensity $1$, independent on $\xi$. We call the random process $N = \{N_1(\xi \lambda(t)), t \geq 0\}$ a {\bf{Mixed Poisson process with Stacy mixing variable}} or {\bf{MPPS-process}}. Briefly $N \in MPPS(\alpha, \beta, c, \lambda(t))$.

\begin{figure}
\begin{minipage}[t]{0.33\linewidth}
   \includegraphics[scale=.55]{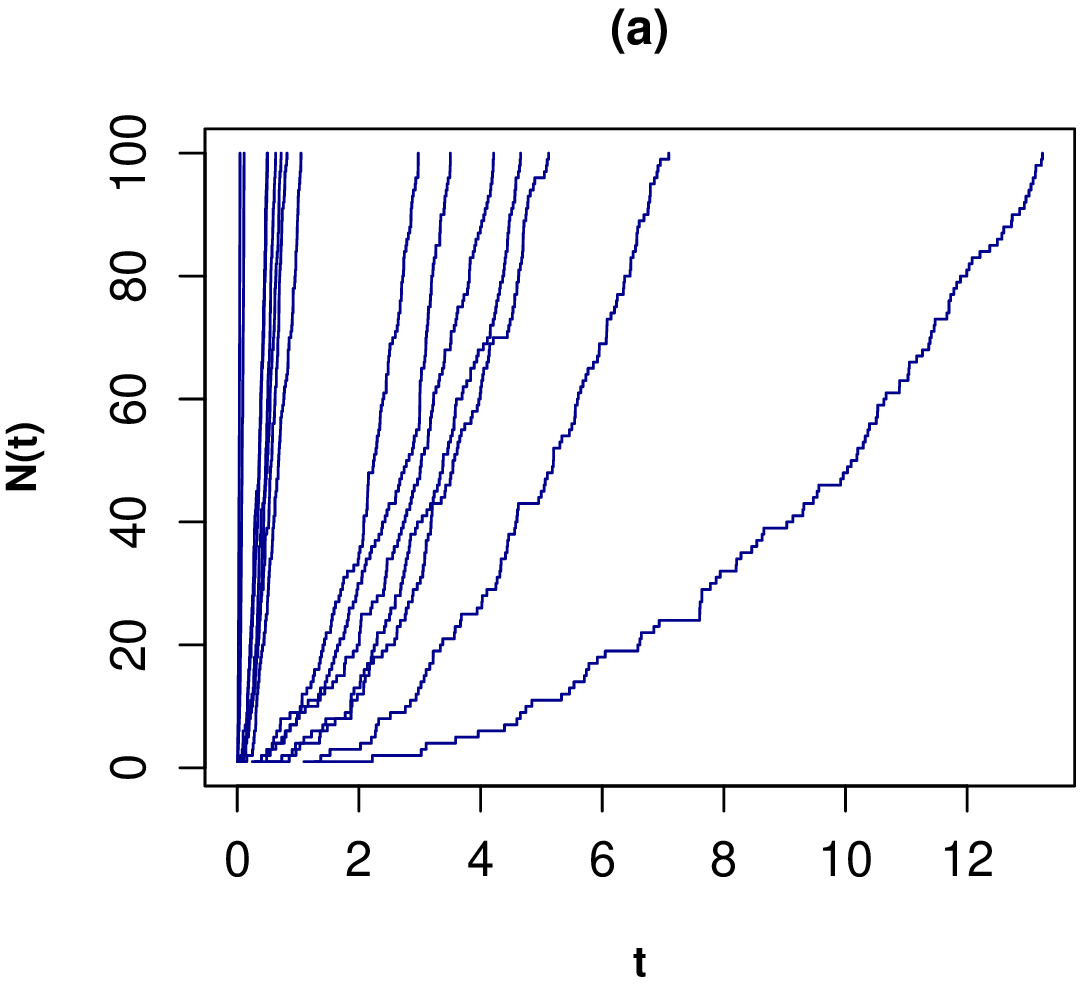}\vspace{-0.3cm}
\end{minipage}
\begin{minipage}[t]{0.33\linewidth}
    \includegraphics[scale=.55]{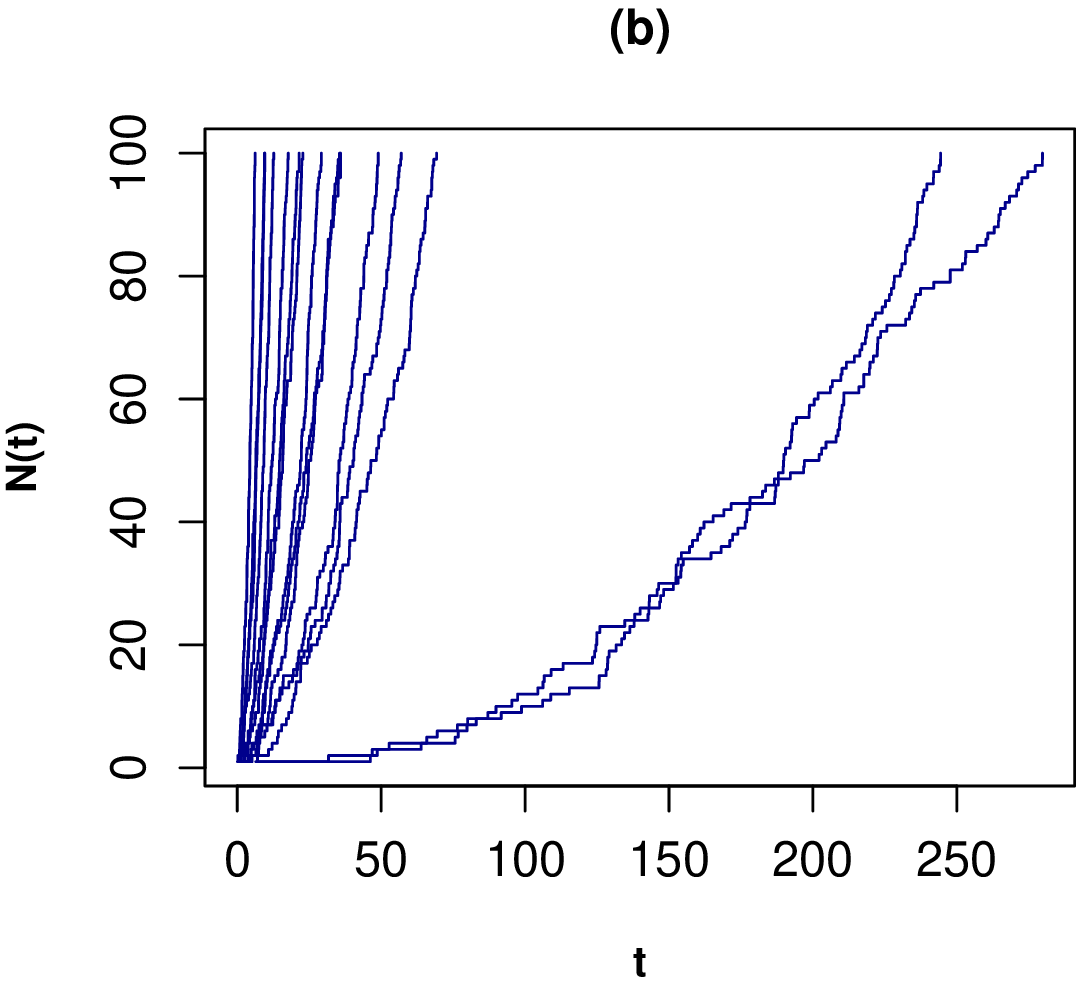}\vspace{-0.3cm}
\end{minipage}
\begin{minipage}[t]{0.33\linewidth}
    \includegraphics[scale=.55]{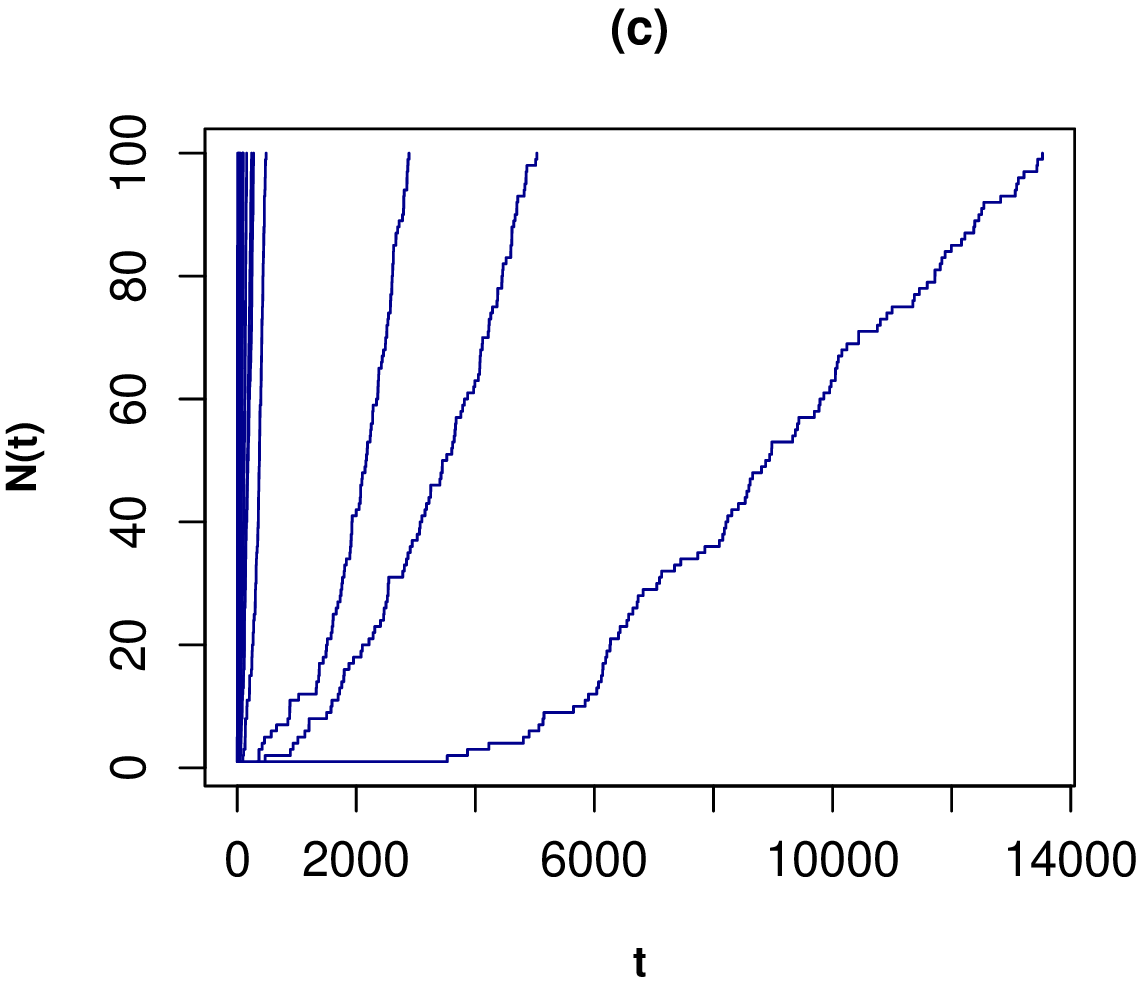}\vspace{-0.3cm}
\end{minipage}
\caption{Fifteen sample paths of $\{N(t), t\geq 0\} \in MPPS(\alpha, \beta, c, t^2)$, for $\alpha = 0.5$, $\beta = 2$, (a)  $c = -0.5$, (b) $c = 1$, (c) $c = 0.5$. \label{fig:pmfMPP1}}
\end{figure}

{\it Notes:} 1. The condition $\lambda(0) = 0$ is not so restrictive as far as otherwise we would shift the Mixed Poisson process with its initial value.

2. When $c = 1$ this is a non-randomly time-changed P$\acute{o}$lya process with its Negative Binomial time intersections
$$\mathbb{P}(N(t) = n) = \mathbb{P}(N_1(\xi \lambda(t)) = n) = \frac{\Gamma(n + \alpha)}{n!\Gamma(\alpha)}\left(\frac{\lambda(t)}{\lambda(t)+\beta}\right)^n \left(\frac{\beta}{\lambda(t)+\beta}\right)^\alpha, \quad n = 0, 1, 2, ..., \quad t > 0. $$

{\bf Definition 7.} \label{Def:7} Let $n \in \mathbb{N}$. We say that a random vector $(N_1, N_2, \ldots, N_n)$ is {\bf Ordered Poisson-Stacy distributed with parameters $\alpha > 0$, $\beta > 0$, $c \in \mathbb{R}\backslash\{0\}$, $0 < \lambda_1 < \lambda_2 < ... < \lambda_n$} if, for all integers $0 \leq k_1 \leq k_2 \leq \ldots \leq k_n$,
$$\mathbb{P}(N_1 = k_1, N_2 = k_2, \ldots, N_n = k_n) = \frac{|c|\lambda_1^{k_1}(\lambda_2 - \lambda_1)^{k_2 - k_1}\ldots(\lambda_n - \lambda_{n-1})^{k_n - k_{n - 1}}}{\beta^{k_n}k_1!(k_2 - k_1)!\ldots(k_n - k_{n-1})!\Gamma(\alpha)}Z_{-c}^{-c\alpha-k_n}\left(\frac{\lambda_n}{\beta}\right),$$
and $\mathbb{P}(N_1 = k_1, N_2 = k_2, \ldots, N_n = k_n) = 0$, otherwise. Briefly, $(N_1, N_2, \ldots, N_n) \in O_{PS}(\alpha, \beta, c; \lambda_1, \lambda_2, ..., \lambda_n)$.

\medskip

{\bf Definition 8.} \label{Def:8} Let $n \in \mathbb{N}$. We say that a random vector $(N_1, N_2, \ldots, N_n)$ is {\bf Mixed Poisson-Stacy distributed with parameters $\alpha > 0$, $\beta > 0$, $c \in \mathbb{R}\backslash\{0\}$, $0 < \lambda_1 < \lambda_2 < ... < \lambda_n$} if, for all $m_1, m_2, \ldots, m_n \in \{0, 1, \ldots\}$,
$$\mathbb{P}(N_1 = m_1, N_2 = m_2, \ldots, N_n = m_n) = \frac{|c|\lambda_1^{m_1}(\lambda_2 - \lambda_1)^{m_2}\ldots(\lambda_n - \lambda_{n-1})^{m_n}}{\beta^{m_1 + m_2 + \ldots + m_n}m_1!m_2!\ldots m_n!\Gamma(\alpha)}Z_{-c}^{-c\alpha- m_1 - m_2 - \ldots - m_n}\left(\frac{\lambda_n}{\beta}\right),$$
and $\mathbb{P}(N_1 = m_1, N_2 = m_2, \ldots, N_n = m_n) = 0$, otherwise. Briefly, $(N_1, N_2, \ldots, N_n) \in M_{PS}(\alpha, \beta, c; \lambda_1, \lambda_2, ..., \lambda_n)$.

\medskip

{\bf Proposition 1.} If $(N_1, N_2, \ldots, N_n) \in O_{PS}(\alpha, \beta, c; \lambda_1, \lambda_2, ..., \lambda_n)$, then
$$(N_1, N_2 - N_1, \ldots, N_n - N_{n-1}) \in M_{PS}(\alpha, \beta, c; \lambda_1, \lambda_2, ..., \lambda_n).$$

{\bf Proposition 2.}  If $(N_1, N_2, \ldots, N_n) \in M_{PS}(\alpha, \beta, c; \lambda_1, \lambda_2, ..., \lambda_n)$, then
$$(N_1, N_1 + N_2, \ldots, N_1 + N_2 + \ldots + N_n) \in O_{PS}(\alpha, \beta, c; \lambda_1, \lambda_2, ..., \lambda_n).$$

Let us consider MPPS-process $\{N_{\alpha, \beta, c, \lambda(t)}, t\geq 0\}$ as a counting process. In the next theorem we denote by $\tau_1, \tau_2, \ldots$ the inter-occurrence (inter-arrival) times, and for $n \in \mathbb{N}$ we denote by $T_n$ the moment of occurrence (arrival) of the $n$-th event. More precisely,
$$T_n = \tau_1 + \tau_2 + \ldots + \tau_n, \quad {\text{and}} \quad \{N(t), \quad t \geq 0\} = \{sup\{i \geq 0: t_i \leq t\}, t \geq 0\}.$$

{\bf Theorem 3.} Let $\lambda(t): [0, \infty) \to [0, \infty)$ be a nonnegative, strictly increasing and continuous function, and $\{N(t), t \geq 0\} \in MPPS(\alpha, \beta, c, \lambda(t))$. Then,

\begin{description}
\item[a)] for all $t > 0$, $N(t) \in MPStacy(\alpha, \frac{\beta}{\lambda(t)}, c)$.
\item[b)] The mean and the variance of this process are
$$\mathbb{E}N(t) = \frac{\Gamma\left(\alpha + \frac{1}{c}\right)}{\beta\Gamma(\alpha)}\lambda(t), \quad \alpha > -\frac{1}{c}$$
$$\mathbb{D}N(t) = \lambda(t)\frac{\Gamma\left(\alpha + \frac{1}{c}\right)}{\beta\Gamma(\alpha)}\left(1 + \frac{\Gamma\left(\alpha + \frac{2}{c}\right)}{\Gamma\left(\alpha + \frac{1}{c}\right)} - \frac{\Gamma\left(\alpha + \frac{1}{c}\right)}{\beta\Gamma(\alpha)}\right), \quad \alpha > -\min\left(\frac{1}{c},\frac{2}{c}\right).$$
\item[c)] The probability generating function of the time intersections is
$$\mathbb{E}(z^{N(t)}) = \frac{|c|}{\Gamma(\alpha)} Z_{-c}^{-c\alpha}\left(\frac{\lambda(t)}{\beta}(1-z)\right), \quad z < 1.$$
\item[d)] For $t > 0$, and $n = 0, 1, \ldots$,
$$P_{\xi}(x|N(t) = n) = \frac{\beta^{c\alpha+n}x^{c\alpha+n-1}}{Z_{-c}^{-c\alpha-n}\left(\frac{\lambda(t)}{\beta}\right)}e^{-\lambda(t)x-(\beta x)^c}, \quad x > 0.$$
\item[e)]  For $t > 0$, and $n = 0, 1, \ldots$, the mean square regression
$$\mathbb{E}(\xi|N(t) = n) = \frac{Z_{-c}^{-c\alpha-n - 1}\left(\frac{\lambda(t)}{\beta}\right)}{\beta Z_{-c}^{-c\alpha-n}\left(\frac{\lambda(t)}{\beta}\right)}.$$
\item[f)] For all $k = 0, 1, \ldots$, $\frac{k}{c}>-\alpha$,
$$\mathbb{E}[N(t)(N(t)-1)(N(t) - k + 1)] = [\lambda(t)]^k\mathbb{E}(\xi^k) = \frac{[\lambda(t)]^k\Gamma\left(\alpha + \frac{k}{c}\right)}{\beta^k\Gamma(\alpha)}.$$
\item[g)] For all $n \in \mathbb{N}$, and $0 \leq t_1 \leq t_2 \leq \ldots \leq t_n$, the f.d.ds. of the random process $\{N(t), t\geq 0\} \in MPPS(\alpha, \beta, c, \lambda(t))$ are
$(N(t_1), N(t_2), \ldots, N(t_n)) \in O_{PS}(\alpha, \beta, c; \lambda(t_1), \lambda(t_2), ..., \lambda(t_n))$.
\item[h)] For all $n \in \mathbb{N}$, and $0 \leq t_1 \leq t_2 \leq \ldots \leq t_n$, the joint distribution of the additive increments of the random process $\{N(t), t\geq 0\} \in MPPS(\alpha, \beta, c, \lambda(t))$ are
$$(N(t_1), N(t_2) - N(t_1), \ldots, N(t_n) - N(t_{n-1})) \in M_{PS}(\alpha, \beta, c; \lambda(t_1), \lambda(t_2), ..., \lambda(t_n)).$$
\item[i)] $\tau_1, \tau_2, \ldots$ are dependent and $Exp-St(\alpha, \beta; c)$ distributed;
\item[j)] $T_n \in Erlang-St(n;\alpha, \beta; c)$;
\end{description}

{\bf Proof:} a) Consider $t > 0$ and $n \in \{0, 1, \ldots\}$. By Definition 5, the independence between the random process $N_1$ and the r.v. $\xi$,  after making substitution $x\beta = y$, and finally using Definition 6 we obtain,
\begin{eqnarray*}
  \mathbb{P}(N(t) = n) &=& \mathbb{P}(N_1(\xi \lambda(t))= n) = \int_0^\infty \mathbb{P}(N_1(\xi \lambda(t))= n|\xi=x)P_\xi(x) dx  \\
   &=& |c|\frac{\beta^{c\alpha}}{\Gamma(\alpha)} \int_0^\infty \mathbb{P}(N_1(x\lambda(t))= n)x^{c\alpha-1}e^{-(\beta x)^c}dx  = |c|\frac{\beta^{c\alpha}}{\Gamma(\alpha)} \int_0^\infty \frac{(x\lambda(t))^n}{n!}e^{-x\lambda(t)}x^{c\alpha-1}e^{-(\beta x)^c}dx\\
   &=& |c|\frac{\beta^{c\alpha}(\lambda(t))^n}{n!\Gamma(\alpha)} \int_0^\infty x^{c\alpha + n -1}e^{-(\beta x)^c-x\lambda(t)}dx = |c|\frac{(\lambda(t))^n}{n!\beta^n\Gamma(\alpha)} \int_0^\infty y^{c\alpha + n -1}e^{-y^c-\lambda(t)\frac{y}{\beta}}dy \\
   &=& |c|\frac{(\lambda(t))^n}{n!\beta^n\Gamma(\alpha)} Z_{-c}^{-c\alpha - n} \left(\frac{\lambda(t)}{\beta}\right).
\end{eqnarray*}

b) All Mixed Poisson distribution is overdispersed. The general formulae for the mean and the variance of Mixed Poisson distribution p. 3, Grandel \cite{GrandelMixed}, and (\ref{moments}) lead us to the desired result.

c) Let $z < 1$ and $t \geq 0$.

\begin{eqnarray*}
 \mathbb{E}(z^{N(t)}) &=& \mathbb{E}(z^{N_1(\xi \lambda(t))}) = \int_0^\infty  \mathbb{E}(z^{N_1(x \lambda(t))}|\xi=x)P_\xi(x) dx = |c|\frac{\beta^{c\alpha}}{\Gamma(\alpha)} \int_0^\infty\mathbb{E}(z^{N_1(x \lambda(t))})x^{c\alpha-1}e^{-(\beta x)^c}dx\\
   &=&  |c|\frac{\beta^{c\alpha}}{\Gamma(\alpha)} \int_0^\infty e^{-x\lambda(t)(1-z)}x^{c\alpha-1}e^{-(\beta x)^c}dx\\
   &=& \frac{|c|}{\Gamma(\alpha)} \int_0^\infty y^{c\alpha-1} e^{-y\frac{\lambda(t)}{\beta}(1-z)-y^c}dy = \frac{|c|}{\Gamma(\alpha)} Z_{-c}^{-c\alpha}\left(\frac{\lambda(t)}{\beta}(1-z)\right).
\end{eqnarray*}

d) The proof follows by the Bayes rule, Definition 6, a), (\ref{pdfStacy}), and (\ref{MixedPoissonStacy}).

f) follows by Remark 2.1, p. 15 in Grandel \cite{GrandelMixed}, a), and (\ref{moments}).

g) and h) are corollaries of the Total probability formula, Definition 6, a), and the independent and Poisson distributed increments of the Poisson process.

i) is a corollary of Definition 6, the fact that the inter-occurrence times is a Homogeneous Poisson process are exponential and Theorem 1, a).

j) In order to prove this statement we use Definition 6, the fact that the occurrence times is a Homogeneous Poisson process are Erlang distributed and Theorem 2, d).

 \hfill $\Box$

{\it{Notes:}} 1. Due to the scale property of the Stacy distribution we can consider $\xi^*(t) \in Stacy(\alpha, \frac{\beta}{\lambda(t)}; c)$, and then, for all $t > 0$,
$$N(t) \stackrel{d}{=} N_1(\xi \lambda(t)) = N_1(\xi^*(t)).$$

2. For $\xi \in Stacy(\alpha, \beta; c)$, by using Definition 6, a), c) and (\ref{StacyMGF}) we obtain the well-known property of Mixed Poisson processes that
$$\mathbb{E}z^{N(t)} = \mathbb{E} e^{-\xi \lambda(t)(1-z)}, \quad z < 1.$$

{\bf Theorem 4.} Suppose that $\alpha > 0$, $\beta > 0$, $c \in \mathbb{R}\backslash\{0\}$, $\xi \in Stacy(\alpha, \beta; c)$, and given $\xi = x$, $x > 0$, the random processes $N_{1, \xi}, N_{2, \xi}, \ldots$ are independent homogeneous Poisson processes with intensity $x$, then
\begin{description}
\item[a)] for all $s \in \mathbb{N}$,
$$N_{1, \xi}(\lambda(t)) + N_{2, \xi}(\lambda(t)) + \ldots + N_{s, \xi}(\lambda(t)) \stackrel{f.d.d.}{=} N_{\alpha, \frac{\beta}{s}, c, \lambda(t)},$$
where $\{N_{\alpha, \frac{\beta}{s}, c, \lambda(t)}, t\geq 0\} \in MPPS(\alpha, \frac{\beta}{s}, c, \lambda(t)).$
\item[b)] If additionally $\eta$ is a discrete r.v. with possible values $1, 2, \ldots$, then
$$N_{1, \xi}(\lambda(t)) + N_{2, \xi}(\lambda(t)) + \ldots + N_{\eta, \xi}(\lambda(t)) \stackrel{f.d.d.}{=} N_{\alpha, \frac{\beta}{\eta}, c, \lambda(t)},$$
where given $\eta = s$, the random process $\{N_{\alpha, \frac{\beta}{\eta}, c, \lambda(t)}, t\geq 0\} \in MPPS(\alpha, \frac{\beta}{s}, c, \lambda(t)).$
\end{description}

{\bf Proof:} a) For $0 < t_1 < t_2$ we are going to consider only bivariate distributions. In an analogous way we can prove the equality of any finite dimensional distributions (f.d.ds.) of these processes. By consequtive application of the definition of a probability generating function (p.g.f.), the formula for the double expectation, the monotonicity of the function $\lambda$, the independent increments of the considered homogeneous Poisson processes and further on the same facts, however in the opposite order:
\begin{eqnarray*}
  & \mathbb{E}&\left[z_1^{N_{1, \xi}(\lambda(t_1)) + N_{2, \xi}(\lambda(t_1)) + \ldots + N_{s, \xi}(\lambda(t_1))}z_2^{N_{1, \xi}(\lambda(t_2)) + N_{2, \xi}(\lambda(t_2)) + \ldots + N_{s, \xi}(\lambda(t_2))-N_{1, \xi}(\lambda(t_1)) - N_{2, \xi}(\lambda(t_1)) - \ldots - N_{s, \xi}(\lambda(t_1))}\right] \\
   &=& \int_0^\infty \mathbb{E}\left[z_1^{N_{1, x}(\lambda(t_1)) + N_{2, x}(\lambda(t_1)) + \ldots + N_{s, x}(\lambda(t_1))}z_2^{N_{1, x}(\lambda(t_2)) + N_{2, x}(\lambda(t_2)) + \ldots + N_{s, x}(\lambda(t_2))-N_{1, x}(\lambda(t_1)) - N_{2, x}(\lambda(t_1)) - \ldots - N_{s, x}(\lambda(t_1))}|\xi=x\right] P_\xi(x) dx\\
   &=& \int_0^\infty \mathbb{E}\left[z_1^{N_{1, x}(\lambda(t_1)) + N_{2, x}(\lambda(t_1)) + \ldots + N_{s, x}(\lambda(t_1))}z_2^{N_{1, x}(\lambda(t_2)) + N_{2, x}(\lambda(t_2)) + \ldots + N_{s, x}(\lambda(t_2))-N_{1, x}(\lambda(t_1)) - N_{2, x}(\lambda(t_1)) - \ldots - N_{s, x}(\lambda(t_1))}\right] P_\xi(x) dx\\
   &=& \int_0^\infty \mathbb{E}\left[z_1^{N_{1, x}(\lambda(t_1)) + N_{2, x}(\lambda(t_1)) + \ldots + N_{s, x}(\lambda(t_1))}\right]  \mathbb{E}\left[z_2^{N_{1, x}(\lambda(t_2)) + N_{2, x}(\lambda(t_2)) + \ldots + N_{s, x}(\lambda(t_2))-N_{1, x}(\lambda(t_1)) - N_{2, x}(\lambda(t_1)) - \ldots - N_{s, x}(\lambda(t_1))}\right] P_\xi(x) dx\\
   &=& \int_0^\infty \left[\mathbb{E}\left[z_1^{N_{1, x}(\lambda(t_1))}\right] \mathbb{E}\left[z_2^{N_{1, x}(\lambda(t_2)) - N_{1, x}(\lambda(t_1))}\right]\right]^s  P_\xi(x) dx = \int_0^\infty \left[\mathbb{E}\left[z_1^{N_{1, x}(\lambda(t_1))} z_2^{N_{1, x}(\lambda(t_2)) - N_{1, x}(\lambda(t_1))}\right]\right]^s  P_\xi(x) dx\\
   &=& \int_0^\infty e^{-sx\lambda(t_1)(1-z_1)} e^{-sx(\lambda(t_2) - \lambda(t_1))(1-z_2)}  P_\xi(x) dx = \int_0^\infty e^{-sx[\lambda(t_1)(1-z_1)+(\lambda(t_2) - \lambda(t_1))(1-z_2)]} P_\xi(x) dx\\
   &=& \int_0^\infty e^{-y[\lambda(t_1)(1-z_1)-(\lambda(t_2) - \lambda(t_1))(1-z_2)]} P_\xi\left(\frac{y}{s}\right) d\frac{y}{s} = \int_0^\infty \mathbb{E}\left[z_1^{N_{1, y}(\lambda(t_1))}\right]\mathbb{E}\left[z_2^{N_{1, y}(\lambda(t_2)) - N_{1, y}(\lambda(t_1))}\right] P_\xi\left(\frac{y}{s}\right) d\frac{y}{s}\\
   &=& \int_0^\infty \mathbb{E}\left[z_1^{N_{1, y}(\lambda(t_1))}z_2^{N_{1, y}(\lambda(t_2)) - N_{1, y}(\lambda(t_1))}\right] P_\xi\left(\frac{y}{s}\right) d\frac{y}{s} = \int_0^\infty \mathbb{E}\left[z_1^{N_{1, y}(\lambda(t_1))}z_2^{N_{1, y}(\lambda(t_2)) - N_{1, y}(\lambda(t_1))}\right] P_{s\xi}(y) dy\\
   &=&  \int_0^\infty \mathbb{E}\left[z_1^{N_{1, y}(\lambda(t_1))}z_2^{N_{1, y}(\lambda(t_2)) - N_{1, y}(\lambda(t_1))}|s\xi = y \right] P_{s\xi}(y) dy =
   \mathbb{E}\left[z_1^{N_{\alpha, \frac{\beta}{s}, c, \lambda(t_1)}}z_2^{N_{\alpha, \frac{\beta}{s}, c, \lambda(t_2)} - N_{\alpha, \frac{\beta}{s}, c, \lambda(t_1)}}\right].
   \end{eqnarray*}

b) The proofs follow by a) and the Total probability formula.

\bigskip
\section{CONCLUSIONS}
Due to the generality of the counting process considered in this work, it is appropriate for modelling of more than just Poisson arrivals. Moreover, as a particular case, it contains the Negative Binomial process. However, an open problem is the one about the estimation of its parameters. We can apply, for example, the method of moments to the inter-arrival times. In that case, we need to invert the Gamma function, and to use this inverse in a system of equations. However, this approach works well only in case when the corresponding moments exist. The conditions for that are described in Theorem 1, d). The situation is similar if we apply, the method of moments to the arrival times. Then, we use Theorem 2, d). If we use the counting process, Theorem 3, b) can help. The system of equations that arises in the maximum likelihood approach has no explicit solution. Therefore, Newton-Like methods have to be used. The explicit analytical form of he estimators is still an open problem.

\section{APPENDIX WITH THE CODES FOR SIMULATION OF THE SAMPLE PATHS}
In the seek of reproducibility of our research, this section contains the code for simulations of the sample paths which were plotted in Figure \ref{fig:pmfMPP1}.

\begin{verbatim}
rm(list=ls(all=TRUE))
n = 100; alpha = 1.5; beta = 0.5; c = 0.5
stacy = rgamma(15, alpha, rate=beta^c)^(1/c);plot(stacy)
exp = mat.or.vec(15, n);
TL = mat.or.vec(15, n);
T = mat.or.vec(15, n);
for (j in 1:15)
    exp[j,]=-log(1-runif(n))/stacy[j]
for (j in 1:15)
  for (i in 1:n)
    TL[j,i]<-sum(exp[j,1:i])
for (j in 1:15)
  for (i in 1:n)
    T[j,1:n] = TL[j,1:n]^(1/2)

m=max(T)

plot(T[1,1:n],1:100, type = "s", xlim=c(0,m), cex.axis=1.2, col = "white",
+ xlab=substitute(paste(bold("t"))),ylab=substitute(paste(bold("N(t)"))),main="(c)")
for (i in 1:15)
lines(T[i,1:n],1:n, type = "s", col = "darkblue")
\end{verbatim}

It allows us to use the Monte-Carlo method and to make an empirical check of the results.

\bigskip
\section{ACKNOWLEDGMENTS} The work was supported by the Scientific Research Fund in Konstantin Preslavsky University of Shumen, Bulgaria under Grant  RD-08-35/18.01.2023; and ANID Chile – CONCYTEC Peru under Grant COVBIO0003.

\bibliographystyle{aipnum-cp}%

\end{document}